\documentclass[11pt]{amsart}

\usepackage{geometry}
\usepackage{todonotes}
\usepackage{amsfonts, amssymb, amscd}
\numberwithin{equation}{section}

\usepackage{bm}
\usepackage{verbatim}
\usepackage{mathrsfs}
\usepackage{graphicx}
\usepackage{tikz-cd}
\usepackage{subcaption}
\usepackage{listings}
\usepackage{subfiles}
\usepackage[toc,page]{appendix}
\usepackage{mathtools}
\usepackage{comment}
\usepackage{enumerate}
\usepackage{enumitem}
\usepackage[all]{xy}

\usepackage{graphicx}
\graphicspath{{images/}}

\usepackage{appendix}
\usepackage{hyperref}
\hypersetup{
    colorlinks=true,
    citecolor=red,
    linkcolor=blue,
    filecolor=magenta,      
    urlcolor=red,
}
\newcommand{\bir}{\dashrightarrow}

\newcommand{\bb}{\bm{b}}
\newcommand{\Mm}{{\bf{M}}}

\newcommand{\Dd}{{\bf{D}}}

\newcommand{\Cc}{\mathbb{C}}

\newcommand{\Qq}{\mathbb{Q}}

\newcommand{\Rr}{\mathbb{R}}

\newcommand{\Span}{\operatorname{Span}}

\newcommand{\Center}{\operatorname{center}}

\newcommand{\Exc}{\operatorname{Exc}}

\newcommand{\rk}{\operatorname{rank}}

\newcommand{\lct}{\operatorname{lct}}

\newcommand{\fol}{\operatorname{fol}}

\newcommand{\Supp}{\operatorname{Supp}}

\newcommand{\mult}{\operatorname{mult}}

\newcommand{\Ff}{\mathcal{F}}

\newcommand{\Ii}{\Gamma}

\newcommand{\Ee}{\mathcal{E}}

\newcommand{\Sing}{\mathrm{Sing}}

\newtheorem{thm}{Theorem}[section]

\newtheorem{cor}[thm]{Corollary}
\newtheorem{lem}[thm]{Lemma}
\newtheorem{prop}[thm]{Proposition}

\newtheorem{claim}[thm]{Claim}

\theoremstyle{definition}
\newtheorem{defn}[thm]{Definition}
\newtheorem{ques}[thm]{Question}
\theoremstyle{definition}

\newtheorem{ex}[thm]{Example}
\newtheorem{nota}[thm]{Notation}

\theoremstyle{definition}

\begin{document}

\title{ACC for lc thresholds for algebraically integrable foliations}
\author{Omprokash Das, Jihao Liu, and Roktim Mascharak}

\date{\today}

\begin{abstract}
We prove the ACC for lc thresholds and the global ACC for algebraically integrable foliations and provide applications.
\end{abstract}

\address{School of Mathematics, Tata Institute of Fundamental Research, Homi Bhabha Road, Navy Nagar, Colaba, Mumbai 400005, India}
\email{omdas@math.tifr.res.in}

\address{Department of Mathematics, Northwestern University, 2033 Sheridan Road, Evanston, IL 60208, USA}
\email{jliu@northwestern.edu}

\address{School of Mathematics, Tata Institute of Fundamental Research, Homi Bhabha Road, Navy Nagar, Colaba, Mumbai 400005, India}
\email{mascharakroktim99@gmail.com}

\maketitle

\tableofcontents

\section{Introduction}\label{sec:Introduction}

We work over the field of complex numbers $\mathbb{C}$. 

Foliations are a well-studied object in algebraic geometry, exhibiting favorable properties in birational geometry and playing a critical role in the minimal model program, most notably in Miyaoka's proof of the abundance conjecture in dimension three \cite{Miy87} (also see \cite[Chapter 9]{Kol+92}). The theory of foliations has been more crucially applied in the abundance conjecture for K\"ahler threefolds  by the first author and Ou \cite{DO23a,DO23b} very recently (see \cite{CHP16,CHP23} for a previous conjectural approach).

In recent years, it has been discovered that many results in classical birational geometry can be extended to foliations, especially in low dimensions. The foundations of the minimal model program have been established for foliated surfaces (cf. \cite{McQ08,Bru15}) and foliated threefolds (cf. \cite{CS20,Spi20,CS21,SS22}). Moreover, many classical questions from the minimal model programs, such as the ascending chain condition (ACC) conjecture for minimal log discrepancies and the ACC conjecture for lc thresholds, have been generalized to foliations and proved in dimensions $2$ and $3$ \cite{Che22,Che23,LLM23,LMX23a,LMX23b}.

With the foundations of the minimal model program for foliations settled in dimensions $\leq 3$, it becomes very interesting to study the minimal model program for foliations in high dimensions, i.e., in dimensions $\geq 4$. However, the foundation of the minimal model program for foliations in high dimensions, particularly the cone theorem, the contraction theorem, and the existence of flips, all seem to be very difficult questions and are currently out of reach using existing methods. A major obstacle is the unknown existence of a log resolution for foliations in dimension $\geq 4$. Nevertheless, in recent years, one particular class of foliations, namely algebraically integrable foliations, has come to our attention. Algebraically integrable foliations are foliations induced by dominant maps. Many foliations satisfying nice numerical properties are known to be algebraically integrable, and it has become one major type of foliation that has been extensively studied, cf. \cite{AD13,AD14,Dru17,AD19,Dru21,DO22,Liu23}.   

Recently, \cite{ACSS21} proved the cone theorem for algebraically integrable foliations \cite[Theorem 3.9]{ACSS21}, and introduced the concept of ``Property $(*)$ modifications", an analogue of dlt modifications for algebraically integrable foliations. Later, in \cite{CS23}, the authors proved the contraction theorem and the existence of flips for algebraically integrable foliations provided that all sequences of klt flips terminate. These results indicate that, from the point of view of the minimal model program, the behavior of algebraically integrable foliations is much closer to that of usual varieties than that of arbitrary foliations. In particular, the structure of algebraically integrable foliations is known to have a tight connection with the canonical bundle formula \cite{ACSS21}, and is essentially used in the proof of the canonical bundle formula over curves (modulo the minimal model program) in characteristic $p>0$ \cite{Ben23}. The foliations structures used in the proof of the abundance conjecture for K\"ahler threefolds \cite{DO23a,DO23b} are also algebraically integrable.

In this paper, we study an interesting question for algebraically integrable foliations: the ascending chain condition (ACC) for lc thresholds.

\begin{thm}[ACC for lc thresholds for algebraically integrable foliations]\label{thm: acc lct alg int foliation}
Let $r$ be a positive integer and $\Ii\subset [0,1]$ a DCC set. Then there exists an ACC set $\Ii'$ depending only on $r$ and $\Ii$ satisfying the following.

Let $(X,\Ff,B)$ be an lc foliated triple of rank $r$ such that $\Ff$ is algebraically integrable, $B\in\Ii$ (i.e. the coefficients of $B$ are contained in $\Ii$), and $D\geq 0$ and $\Rr$-Cartier $\Rr$-divisor such that $D\in\Ii$. Then the lc threshold
$$\lct(X,\Ff,B;D):=\sup\{t\mid t\geq 0, (X,\Ff,B+D)\text{ is lc}\}$$
belongs to $\Ii'$.
\end{thm}
We note that the set $\Ii'$ depends on $r$ and $\Ii$ rather than $\dim X$ and $\Ii$. Furthermore, for foliations in dimensions $\leq 3$, Theorem \ref{thm: acc lct alg int foliation} has been proven in \cite{Che22}, even for those that are not algebraically integrable. 

For usual pairs, where $\Ff=T_X$, Theorem \ref{thm: acc lct alg int foliation} has been established in \cite[Theorem 1.1]{HMX14}. This is a fundamental result in birational geometry, as it implies the DCC of volumes and effective birationality \cite[Theorem 1.3]{ACSS21}, and has numerous applications in the minimal model program, such as the termination of effective flips \cite{Bir07} and the existence of lc flips \cite{Bir12}. Hence, we anticipate that the ACC for lc thresholds in algebraically integrable foliations will also be valuable for the minimal model program of foliations.

As a companion result of  Theorem \ref{thm: acc lct alg int foliation}, we prove the following theorem, known as the global ACC:

\begin{thm}[Global ACC for algebraically integrable foliations]\label{thm: global acc alg int foliation}
Let $r$ be a positive integer and $\Ii\subset [0,1]$ a DCC set. Then there exists a finite set $\Ii_0\subset\Ii$ depending only on $r$ and $\Ii$ satisfying the following.

Let $(X,\Ff,B)$ be an lc foliated triple of rank $r$ such that $\Ff$ is algebraically integrable, $B\in\Ii$, and $K_{\Ff}+B\equiv 0$. Then $B\in\Ii_0$.
\end{thm}
We note once again that the set $\Ii_0$ depends on $r$ and $\Ii$ rather than $\dim X$ and $\Ii$. Additionally, for foliations (possibly non-algebraically integrable), Theorem \ref{thm: global acc alg int foliation} has been proven in \cite{Che22} for $\dim X=2$ and in \cite{LLM23,LMX23b} for $\dim X=3$.

For usual pairs where $\Ff=T_X$, Theorem \ref{thm: global acc alg int foliation} has been established in \cite[Theorem 1.5]{HMX14}. This result has numerous applications, especially in the study of boundedness results for varieties admitting log Calabi-Yau structures (cf. \cite{HX15,Bir19,Bir21}).

As a straightforward corollary of Theorem \ref{thm: global acc alg int foliation}, we obtain the following result:

\begin{cor}[Global ACC for rank one foliations]\label{cor: global acc rank 1 foliation}
Let $\Ii\subset [0,1]$ be a DCC set. Then there exists a finite set $\Ii_0\subset\Ii$ depending only on $\Ii$ satisfying the following.

Let $(X,\Ff,B)$ be an lc foliated triple of rank $1$ such that $B\in\Ii$ and $K_{\Ff}+B\equiv 0$. Then $B\in\Ii_0$.
\end{cor}

In addition to our main theorems, we also obtain a result on the (log) abundance of algebraically integrable foliations. It is worth noting that \cite[Proposition 4.24]{Dru21} has proven Theorem \ref{thm: log abundance algebraically integrable foliation} for the special case when $\Ff$ has canonical singularities and $B=0$.

\begin{thm}\label{thm: log abundance algebraically integrable foliation} 
Let $(X,\Ff,B)$ be an lc foliated triple such that $\Ff$ is algebraically integrable, $B\in\Ii$, and $K_{\Ff}+B\equiv 0$. Then $K_{\Ff}+B\sim_{\mathbb R}0$.
\end{thm}

It is also worth noting that Theorem \ref{thm: log abundance algebraically integrable foliation} has been proven for foliations (possibly non-algebraically integrable) in dimension $\leq 3$ in \cite{LLM23} (\cite{CS20,CS21} for the case of $\Qq$-divisors).

As a direct consequence of the ACC for lc thresholds in algebraically integrable foliations, we establish the existence of uniform lc rational polytopes for such foliations.

\begin{thm}\label{thm: uniform rational polytope foliation intro}
Let $r$ be a positive integer, $v_1^0,\dots,v_m^0$ positive integers, and $\bm{v}_0:=(v_1^0,\dots,v_m^0)$. Then there exists an open set $U\ni \bm{v}_0$ of the rational envelope of $\bm{v}_0$ depending only on $r$ and $\bm{v}_0$ satisfying the following. 

Let $(X,\Ff,B=\sum_{j=1}^mv_j^0B_j)$ be a foliated lc triple of rank $r$, such that $\Ff$ is algebraically integrable, $\rk\Ff=r$, and $B_j\geq 0$ are distinct Weil divisors. Then  $(X,\Ff,B=\sum_{j=1}^mv_jB_j)$ is lc for any $(v_1,\dots,v_m)\in U$.
\end{thm}
We note that Theorem \ref{thm: uniform rational polytope foliation intro} has been proven for foliations (possibly non-algebraically integrable) in dimension $\leq 3$ in \cite{LMX23a,LMX23b}, while for usual pairs (when $\Ff=T_X$), it has been established in \cite{HLS19}. Despite their technicality, uniform lc rational polytopes are powerful tools in birational geometry and have numerous applications, such as the ACC conjecture for minimal log discrepancies and the boundedness of complements. They are also key ingredients in proving the global ACC for foliated threefolds \cite{LMX23b}.

As an application of Theorem \ref{thm: uniform rational polytope foliation intro}, we obtain a result on the accumulation points of lc thresholds for algebraically integrable foliations:

\begin{cor}\label{cor: accumulation point ai foliation lct}
Let $r$ be a positive integer and $\Ii\subset [0,1]$ be a DCC set such that $\bar\Ii\subset\mathbb Q$. We consider the set of foliated lc thresholds
$$\lct_{a.i.\fol,}(r,\Ii):=\{\lct(X,\Ff,B;D)\mid \rk\Ff=r,\Ff\text{ is algebraically integrable}, B\in\Ii, D\in\mathbb N^+\}.$$
Then the accumulation points of $\lct_{\fol}(r,\Ii)$ are rational numbers.
\end{cor}

Moreover, we anticipate obtaining a more detailed characterization of the accumulation points for lc thresholds. For instance, a natural question to ask is the following:
\begin{ques}\label{ques: accumulation point of foliation lct}
Suppose that $1\in\Ii$, $1$ is the only accumulation point of $\Ii$, and $\Ii=\Ii_+\cap [0,1]$.
\begin{enumerate}
    \item Are the accumulation points of $\lct_{a.i.\fol,}(r,\Ii)$ equal to $\lct_{a.i.\fol,}(r,\Ii)\backslash\{1\}$?
    \item Is $\lct_{a.i.\fol,}(r,\Ii)$ standardized, in the sense of \cite{LMX22}?
\end{enumerate}
\end{ques}
For the case of usual pairs where $\Ff=T_X$, Question \ref{ques: accumulation point of foliation lct}(1) has been established in \cite[Theorem 1.11]{HMX14}, while Question \ref{ques: accumulation point of foliation lct}(2) has been addressed in \cite[Theorem 4.2]{LMX22}. 

P. Cascini asked us the following more aggressive question:

\begin{ques}\label{ques: same lct sets}
Let $r$ be a positive integer, $\Ii\subset [0,1]$ a DCC set, and
$$\lct(d,\Ii):=\{\lct(X,B;D)\mid \dim X=d,B\in\Ii,D\in\mathbb N^+\}.$$
Do we have
$$\lct(r,\Ii)=\lct_{a.i.\fol}(r,\Ii)$$
when $\Ii$ is the standard set $\{1-\frac{1}{n}\mid n\in\mathbb N^+\}\cup\{1\}$ (or some other DCC set satisfying similar nice properties, e.g. hyperstandard set)?
\end{ques}
It is clear that a positive answer to Question \ref{ques: same lct sets} will imply a positive answer to Question \ref{ques: accumulation point of foliation lct} directly by \cite[Section 11]{HMX14}, at least when $\Ii$ is the standard set. It is also clear that $\lct(r,\Ii)\subset\lct_{a.i.\fol}(r,\Ii)$.

\medskip

\noindent\textit{Structure of the paper}. In Section \ref{sec: Preliminaries}, we provide a concise overview of fundamental concepts and basic results in foliations. In Section \ref{sec: adjunction}, we present explicit adjunction formulas that will be employed later in the paper. In Section \ref{sec: acss modification}, we introduce the concept of ACSS foliated triples and establish the existence of ACSS modifications. These results represent stronger versions of Property $(*)$ foliations and Property $(*)$ modifications as defined in \cite{ACSS21}. In Section \ref{sec: proof of the main theorems}, we prove Theorems \ref{thm: acc lct alg int foliation}, \ref{thm: global acc alg int foliation}, \ref{thm: log abundance algebraically integrable foliation}, and Corollary \ref{cor: global acc rank 1 foliation}. In Section \ref{sec: ulrp}, we prove Theorem \ref{thm: uniform rational polytope foliation intro} and Corollary \ref{cor: accumulation point ai foliation lct}.

\medskip

\noindent\textbf{Acknowledgements}. We thank Paolo Cascini, Guodu Chen, Christopher D. Hacon, Jingjun Han, Yuchen Liu, Fanjun Meng, and Lingyao Xie for useful discussions. We would like to acknowledge the assistance of ChatGPT in polishing the wording.

\section{Preliminaries}\label{sec: Preliminaries}

We will work over the field of complex numbers $\Cc$. Throughout the paper, our main focus will be on normal quasi-projective varieties to ensure consistency with the references. However, it is worth noting that most results should also hold for normal varieties that are not necessarily quasi-projective. Similarly, most of the results in our paper should hold for any algebraically closed field of characteristic zero. We will use the standard notations and definitions as in \cite{KM98, BCHM10} and use them freely. For foliations, we will follow the notations and definitions as in \cite{LLM23}, which are largely consistent with those in \cite{CS20,ACSS21,CS21}. For generalized pairs, we will follow the notations and definitions as in \cite{HL21}.

\subsection{Sets}

\begin{defn}\label{defn: DCC and ACC}
Let $\Ii\subset\Rr$ be a set. We say that $\Ii$ satisfies the \emph{descending chain condition} (DCC) if any decreasing sequence in $\Ii$ stabilizes. We say that $\Ii$ satisfies the \emph{ascending chain condition} (ACC) if any increasing sequence in $\Ii$ stabilizes. We define
$$\Ii_+:=\{0\}\cup\left\{\sum_{i=1}^n \gamma_i\Biggm| n\in\mathbb N^+, \gamma_1,\dots,\gamma_n\in\Ii\right\}.$$
\end{defn}

\begin{defn}
Let $m$ be a positive integer and $\bm{v}\in\mathbb R^m$. The \emph{rational envelope} of $\bm{v}$ is the minimal rational affine subspace of $\mathbb R^m$ which contains $\bm{v}$. For example, if $m=2$ and $\bm{v}=\left(\frac{\sqrt{2}}{2},1-\frac{\sqrt{2}}{2}\right)$, then the rational envelope of $\bm{v}$ is $(x_1+x_2=1)\subset\mathbb R^2_{x_1x_2}$.
\end{defn}

\subsection{Foliations}

\begin{defn}[Foliations, {cf. \cite[Section 2.1]{CS21}}]\label{defn: foliation}
Let $X$ be a normal variety. A \emph{foliation} on $X$ is a coherent sheaf $\Ff\subset T_X$ such that
\begin{enumerate}
    \item $\Ff$ is saturated in $T_X$, i.e. $T_X/\Ff$ is torsion free, and
    \item $\Ff$ is closed under the Lie bracket.
\end{enumerate}
The \emph{rank} of the foliation $\Ff$ is the rank of $\Ff$ as a sheaf and is denoted by $\rk\Ff$. The \emph{co-rank} of $\Ff$ is $\dim X-\rk\Ff$. The \emph{canonical divisor} of $\Ff$ is a divisor $K_\Ff$ such that $\mathcal{O}_X(-K_{\mathcal{F}})\cong\mathrm{det}(\Ff)$. We define $N_{\Ff}:=(T_X/\Ff)^{\vee\vee}$ and $N_{\Ff}^*:=N_{\Ff}^{\vee}$.

If $\Ff=0$, then we say that $\Ff$ is a \emph{foliation by points}.
\end{defn}

\begin{defn}[Singular locus]
     Let $X$ be a normal variety and $\Ff$ a rank $r$ foliation on $X$. We can associate to $\Ff$ a morphism $$\phi: \Omega_X^{[r]}\to \mathcal{O}_X(K_{\Ff})$$ defined by taking the double dual of the $r$-wedge product of the map $\Omega^1_X\to \Ff^{\vee}$, induced by the inclusion $\Ff\to T_X$. This yields a map $$\phi': (\Omega_X^{[r]}\otimes\mathcal{O}_X(-K_{\Ff}))^{\vee\vee}\to \mathcal{O}_X$$ and we define the singular locus, denoted as $\Sing\Ff$, to be the co-support of the image of $\phi'$.
\end{defn}

\begin{defn}[Pullbacks and pushforwards, {cf. \cite[3.1]{ACSS21}}]\label{defn: pullback}
Let $X$ be a normal variety, $\Ff$ a foliation on $X$, $f: Y\dashrightarrow X$ a dominant map, and $g: X\dashrightarrow X'$ a birational map of normal varieties. We denote $f^{-1}\Ff$ the \emph{pullback} of $\Ff$ on $Y$ as constructed in \cite[3.2]{Dru21}. We also say that $f^{-1}\Ff$ is the \emph{induced foliation} of $\Ff$ on $Y$. We define the \emph{pushforward} of $\Ff$ on $X'$ as $(g^{-1})^{-1}\Ff$ and denote it by $g_*\Ff$.
\end{defn}

\begin{defn}[Algebraically integrable foliations, {cf. \cite[3.1]{ACSS21}}]\label{defn: algebraically integrable}
Let $X$ be a normal variety and $\Ff$ a foliation on $X$. We say that $\Ff$ is an \emph{algebraically integrable foliation} if there exists a dominant rational map $f: X\dashrightarrow Y$ such that $\Ff=f^{-1}\Ff_Y$, where $\Ff_Y$ is a foliation by points. In this case, we say that $\Ff$ is \emph{induced by $f$}.
\end{defn}

\begin{lem}\label{lem:foliation-invariant}
    Let $f:X'\to X$ be a proper birational morphism between normal varieties, $\Ff$ a foliation on $X$, and $\Ff':=f^{-1}\Ff$ the pullback foliation of $\Ff$ on $X'$. Then $\Ff'$ is algebraically integrable if and only if $\Ff$ is algebraically integrable.
\end{lem}
\begin{proof}
First we prove the if part. Assume that $\Ff$ is algebraically integrable. Let $\pi: X\dashrightarrow Y$ be the dominat map which induces $\Ff$, then $f\circ\pi: X'\dashrightarrow Y$ is a dominant map which induces $\Ff'$.

Now we prove the only if part. Assume that $\Ff'$ is algebraically integrable. Then there is a dominant rational map $\pi':X'\bir Y'$ which induces $\Ff'$. Let $U'\subset X'$ be the largest open set such that $f|_{U'}: U'\to f(U')\subset X$ is an isomorphism and $\pi'|_{U'}:U'\to Y'$ is a morphism. Then $\Ff|_{f(U')}\cong \Ff'|_{U'}\cong T_{U'/Y'}$. Thus $\Ff$ is induced by the rational map $\pi'\circ f^{-1}: X\bir Y'$.
\end{proof}

The following theorem essentially follows from \cite[Theorem 1.1]{CP19}.
\begin{thm}[{\cite[Theorem 3.1]{LLM23},\cite[Theorem 1.1]{CP19}}]\label{thm: subfoliation algebraic integrable}
Let $\Ff$ be a foliation on a normal projective variety $X$ such that $K_{\Ff}$ is not pseudo-effective. Then there exists an algebraically integrable foliation $\Ee$ such that $0\not=\Ee\subset\Ff$.
\end{thm}

\begin{defn}[Invariant subvarieties, {cf. \cite[3.1]{ACSS21}}]\label{defn: f-invariant}
Let $X$ be a normal variety, $\Ff$ a foliation on $X$, and $S\subset X$ a subvariety. We say that $S$ is $\Ff$-invariant if and only if for any open subset $U\subset X$ and any section $\partial\in H^0(U,\Ff)$, we have $$\partial(\mathcal{I}_{S\cap U})\subset \mathcal{I}_{S\cap U}$$ 
where $\mathcal{I}_{S\cap U}$ is the ideal sheaf of $S\cap U$. 
\end{defn}

\begin{defn}[Tangent and transverse]\label{defn: tangent to foliation}
Let $X$ be a normal variety, $\Ff$ a foliation on $X$, and $V\subset X$ a subvariety. Suppose that $\Ff$ is a foliation induced by a dominant rational map $X\dashrightarrow Z$. 
 We say that $V$ is \emph{tangent} to $\Ff$ if there exists a birational morphism $\mu: X'\rightarrow X$, an equidimensional contraction $f': X'\rightarrow Z$, and a subvariety $V'\subset X'$, such that
    \begin{enumerate}
    \item $\mu^{-1}\Ff$ is induced by $f'$, and
        \item $V'$ is contained in a fiber of $f'$ and $\mu(V')=V$.
    \end{enumerate}
    We say that $V$ is \emph{transverse} to $\Ff$ if $V$ is not tangent to $\Ff$.
\end{defn}

\begin{defn}[Special divisors on foliations, cf. {\cite[Definition 2.2]{CS21}}]\label{defn: special divisors on foliations}
Let $X$ be a normal variety and $\Ff$ a foliation on $X$. For any prime divisor $C$ on $X$, we define $\epsilon_{\Ff}(C):=1$ if $C$ is not $\Ff$-invariant, and  $\epsilon_{\Ff}(C):=0$ if $C$ is $\Ff$-invariant. If $\Ff$ is clear from the context, then we may use $\epsilon(C)$ instead of $\epsilon_{\Ff}(C)$. For any $\Rr$-divisor $D$ on $X$, we define $$D^{\Ff}:=\sum_{C\mid C\text{ is a component of }D}\epsilon_{\Ff}(C)C.$$
Let $E$ be a prime divisor over $X$ and $f: Y\rightarrow X$ a projective birational morphism such that $E$ is on $Y$. We define $\epsilon_{\Ff}(E):=\epsilon_{f^{-1}\Ff}(E)$. It is clear that $\epsilon_{\Ff}(E)$ is independent of the choice of $f$.
\end{defn}

\subsection{Polarized foliations}

\begin{defn}[$\bb$-divisors]\label{defn: b divisors} Let $X$ be a normal quasi-projective variety. We call $Y$ a \emph{birational model} over $X$ if there exists a projective birational morphism $Y\to X$. 

Let $X\dashrightarrow X'$ be a birational map. For any valuation $\nu$ over $X$, we define $\nu_{X'}$ to be the center of $\nu$ on $X'$. A \emph{$\bb$-divisor} $\Dd$ over $X$ is a formal sum $\Dd=\sum_{\nu} r_{\nu}\nu$ where $\nu$ are valuations over $X$ and $r_{\nu}\in\mathbb R$, such that $\nu_X$ is not a divisor except for finitely many $\nu$. If in addition, $r_{\nu}\in\Qq$ for every $\nu$, then $\Dd$ is called a \emph{$\Qq$-$\bb$-divisor}. The \emph{trace} of $\Dd$ on $X'$ is the $\Rr$-divisor
$$\Dd_{X'}:=\sum_{\nu_{X'}\text{ is a divisor}}r_\nu\nu_{X'}.$$
If $\Dd_{X'}$ is $\Rr$-Cartier and $\Dd_{Y}$ is the pullback of $\Dd_{X'}$ on $Y$ for any birational model $Y$ over $X'$, we say that $\Dd$ \emph{descends} to $X'$ and $\Dd$ is the \emph{closure} of $\Dd_{X'}$, and write $\Dd=\overline{\Dd_{X'}}$. 

Let $X\rightarrow U$ be a projective morphism and assume that $\Dd$ is a $\bb$-divisor over $X$ such that $\Dd$ descends to some birational model $Y$ over $X$. If $\Dd_Y$ is a Cartier divisor, then we say that $\Dd$ is \emph{$\bb$-Cartier}. If $\Dd$ can be written as an $\Rr_{\geq 0}$-linear combination of nef$/U$ $\bb$-Cartier $\bb$-divisors, then we say that $\Dd$ is \emph{NQC}$/U$.

We let $\bm{0}$ be the $\bb$-divisor $\bar{0}$.
\end{defn}

\begin{defn}[Generalized foliated quadruples]
A \emph{generalized foliated sub-quadruple} (\emph{sub-gfq} for short) $(X,\Ff,B,\Mm)/U$ consists of a normal quasi-projective variety $X$, a foliation $\Ff$ on $X$, an $\Rr$-divisor $B$ on $X$, a projective morphism $X\rightarrow U$, and an NQC$/U$ $\bb$-divisor $\Mm$ over $X$, such that $K_{\Ff}+B+\Mm_X$ is $\mathbb R$-Cartier. If $B\geq 0$, then we say that $(X,\Ff,B,\Mm)/U$ is a \emph{generalized foliated quadruple} (\emph{gfq} for short). If $U=\{pt\}$, we usually drop $U$ and say that $(X,\Ff,B,\Mm)$ is \emph{projective}. 

Let $(X,\Ff,B,\Mm)/U$ be a (sub-)gfq. If $\Mm=\bm{0}$, then we may denote $(X,\Ff,B,\Mm)/U$ by $(X,\Ff,B)/U$ or $(X,\Ff,B)$, and say that $(X,\Ff,B)$ is a \emph{foliated (sub-)triple} (\emph{f-(sub-)triple} for short). If $\Ff=T_X$, then we may denote $(X,\Ff,B,\Mm)/U$ by $(X,B,\Mm)/U$, and say that $(X,B,\Mm)/U$ is a \emph{generalized (sub-)pair} (\emph{g-(sub-)pair} for short). If $\Mm=\bm{0}$ and $\Ff=T_X$, then we may denote $(X,\Ff,B,\Mm)/U$ by $(X,B)/U$, and say that $(X,B)/U$ is a \emph{(sub-)pair}. 

A (sub-)gfq (resp. f-(sub-)triple, f-(sub-)pair, g-(sub-)pair, (sub-)pair) $(X,\Ff,B,\Mm)/U$ (resp. $(X,\Ff,B)/U$,$(X,B,\Mm)/U$, $(X,B)/U$) is called a \emph{$\mathbb Q$-(sub-)gfq} (resp. \emph{$\mathbb Q$-f-(sub-)triple, $\mathbb Q$-g-(sub-)pair, $\mathbb Q$-(sub-)pair} if $B$ is a $\mathbb Q$-divisor and $\Mm$ is a $\mathbb Q$-$\bb$-divisor.
\end{defn}

\begin{nota}
In the previous definition, if $U$ is not important, we may also drop $U$. This usually happens when we emphasize the structures of $(X,\Ff,B,\Mm)$ which are independent of the choice of $U$, such as the singularities of $(X,\Ff,B,\Mm)$. In addition, if $B=0$, then we may drop $B$.
\end{nota}

\begin{defn}[Singularities of gfqs]\label{defn: gfq singularity}
Let $(X,\Ff,B,\Mm)$ be a (sub-)gfq. For any prime divisor $E$ over $X$, let $f: Y\rightarrow X$ be a birational morphism such that $E$ is on $Y$, and suppose that
$$K_{\Ff_Y}+B_Y+\Mm_Y:=f^*(K_\Ff+B+\Mm_X)$$
where $\Ff_Y:=f^{-1}\Ff$. We define $a(E,\Ff,B,\Mm):=-\mult_EB_Y$ to be the \emph{discrepancy} of $E$ with respect to $(X,\Ff,B,\Mm)$. It is clear that $a(E,\Ff,B,\Mm)$ is independent of the choice of $Y$. If $\Mm=\bm{0}$, then we let $a(E,\Ff,B):=a(E,\Ff,B,\Mm)$. If $\Ff=T_X$, then we let $a(E,X,B,\Mm):=a(E,\Ff,B,\Mm)$. If $\Mm=\bm{0}$ and $\Ff=T_X$, then we let $a(E,X,B):=a(E,\Ff,B,\Mm)$.

We say that $(X,\Ff,B,\Mm)$ is \emph{(sub-)lc} (resp. \emph{(sub-)klt}) if $a(E,\Ff,B,\Mm)\geq -\epsilon_{\Ff}(E)$ (resp. $>-\epsilon_{\Ff}(E)$) for any prime divisor $E$ over $X$. An \emph{lc place} of $(X,\Ff,B,\Mm)$ is a prime divisor $E$ over $X$ such that $a(E,\Ff,B,\Mm)=-\epsilon_{\Ff}(E)$. An \emph{lc center} of $(X,\Ff,B,\Mm)$ is the center of an lc place of $(X,\Ff,B,\Mm)$ on $X.$
\end{defn}

\begin{defn}[Lc threshold]
Let $(X,\Ff,B)$ be an lc f-sub-triple and $D$ an $\Rr$-Cartier $\Rr$-divisor on $X$. An \emph{lc threshold} (\emph{lct} for short) of $D$ with respect to $(X,\Ff,B)$ is a real number $t_0$, such that
\begin{enumerate}
    \item $(X,\Ff,B+t_0D)$ is sub-lc, and
    \item for any positive real number $\delta$, either $(X,\Ff,B+(t_0+\delta)D)$ or  $(X,\Ff,B+(t_0-\delta)D)$ is not sub-lc.
\end{enumerate}
When $D\geq 0$, lc thresholds of $D$ with respect to $(X,\Ff,B)$ are unique, and we denote the lc threshold of $D$ with respect to $(X,\Ff,B)$ by $\lct(X,\Ff,B;D)$. 
\end{defn}

\subsection{Foliated log resolution}
%\footnote{Liu: This subsection was previously at the beginning of Section 4. I move it here because we need it in Proposition \ref{prop: general hyperplane invariant}}

\begin{defn}[{\cite[3.2. Log canonical foliated pairs]{ACSS21}}]\label{defn: foliated log smooth}
Let $(X,\Ff,B)$ be an f-sub-triple such that $\Ff$ is algebraically integrable. We say that $(X,\Ff,B)$ is \emph{foliated log smooth} if there exists a contraction $\pi: X\rightarrow Z$ satisfying the following.
\begin{enumerate}
\item $X$ has at most quotient toric singularities.
\item $\Ff$ is induced by $\pi$.
    \item $(X,\Sigma_X)$ is toroidal for some reduced divisor $\Sigma_X$ such that $\Supp B\subset\Sigma_X$.  In particular, $(X,\Supp B)$ is toroidal, and $X$ is $\Qq$-factorial klt.
    \item There exists a log smooth toroidal pair $(Z,B_Z)$ such that $$\pi: (X,\Supp B)\rightarrow (Z,\Supp B_Z)$$ is an equidimensional  toroidal contraction. 
\end{enumerate}
We say that the contraction $\pi: X\rightarrow Z$ is \emph{associated to} $(X,\Ff,B)$.
\end{defn}

\begin{lem}{\cite[Lemma 3.1]{ACSS21}}\label{lem: foliated log smooth imply lc}
Let $(X,\Ff,B)$ be an f-sub-triple such that $\Ff$ is algebraically integrable and $(X,\Ff,B)$ is foliated log smooth. Then $(X,\Ff,B^\Ff)$ is lc.
\end{lem}

\begin{defn}\label{defn: log resolution}
Let $X$ be a normal variety, $B$ an $\Rr$-divisor on $X$, and $\Ff$ an algebraically integrable foliation on $X$. A \emph{foliated log resolution} of $(X,\Ff,B)$ is a projective birational morphism $f: Y\rightarrow X$ such that $(Y,\Ff_Y:=f^{-1}\Ff,B_Y:=f^{-1}_*B+\Exc(f))$ is foliated log smooth, where $\Exc(f)$ is the reduced exceptional divisor of $f$.
\end{defn}

\section{Precise adjunction formulas for algebraically integrable foliations}\label{sec: adjunction}

In this section we prove some adjunction formulas for algebraically integrable foliations. 

\begin{defn}
Let $X$ be a normal variety, $\Ff$ an algebraically integrable foliation, and $S$ a prime divisor on $X$ with normalization $S^\nu\rightarrow S$. We define the \emph{restricted foliation} $\Ff_S$ on $S$ in the following way:
\begin{enumerate}
    \item If $S$ is not $\Ff$-invariant, then we let $\iota: S^\nu\rightarrow X$ be the natural morphism and $\Ff_S:=\iota^{-1}\Ff$. 
    \item If $S$ is $\Ff$-invariant, then we let $U\subset X$ be an open set which contains the generic point of $S$ and does not intersect $\Sing(\Ff)\cup\Sing(X)\cup\Sing(S)$, and let $S':=S\cap U$. Since $S$ is $\Ff$-invariant, the natural inclusion $T_{\Ff}|_{S'}\rightarrow T_X|_{S'}$ factors through $T_{S'}$, hence defines a foliation on $S'$, which yields a foliation $\Ff_S$ on $S^\nu$. 
\end{enumerate}
\end{defn}
We remark that our definition coincide with the construction in \cite[Proposition 3.2]{ACSS21}.

\subsection{Adjunction preserves algebraic integrability}

The following proposition was mentioned in \cite[Remark 3.4]{ACSS21}. For the reader's convenience, we provide a full proof here.

\begin{prop}\label{prop: a.i preserved adjunction}
Let $\Ff$ be an algebraically integrable foliation, $S$ a prime divisor on $X$, and $\pi:S^\nu\rightarrow S$ the normalization of $S$. Let $\Ff_S$ be the restriction of the foliation $\Ff$ on $S$. Then:
\begin{enumerate}
    \item $\Ff_S$ is algebraically integrable, and
    \item $\rk\Ff_S=\rk\Ff-\epsilon_{\Ff}(S)$.
\end{enumerate}
\end{prop}
\begin{proof}
Let $f:X\dashrightarrow Z$ be the dominant rational map inducing the foliation $\Ff$ on $X$. Suppose that $\overline X$ is the normalization of the graph of $f$, and $p:\overline X\to X$ and $q:\overline X\to Z$ are the projections. Then possibly replacing $X$ by $\overline X$, $\Ff$ by $p^{-1}\Ff$, $S$ by $p^{-1}_*S$, and $f$ by $q$, we may assume that $\Ff$ is induced by a contraction $f: X\to Z$. Let $g:X'\to X$ be a birational morphism of normal varieties and $\Ff'$ is the pullback foliation on $X'$ induced by an  equidimensional contraction $f':X'\to Z'$ as in \cite[Theorem 2.2]{ACSS21}. Let $S':=g^{-1}_*S$, and let $\tau:S'^\nu\to S'$ be the normalization morphism. Let $\theta:S'^\nu\to S^\nu$ be the induced birational morphism. Then $\Ff_{S'}=\theta^{-1}\Ff_S$ and $\theta_*\Ff_{S'}=\Ff_S$.

Now if $S$ is not $\Ff$-invariant, then $S'$ is not $\Ff'$-invariant and thus $f'(S')=Z'$ and $\Ff_{S'}$ is induced by $f'|_{S'}\circ \tau$. Since $\Ff_{S'}=\theta^{-1}\Ff_S$, it follows from Lemma \ref{lem:foliation-invariant} that $\Ff_S$ is algebraically integrable. Since the leaves of $\Ff_S$ are intersection of $S$ and the leaves of $f$, we get (2).

% $f(S)=Z$, and $\Ff_S$ is induced by $f|_S\circ \pi$. So assume that  

% By \cite[Theorem 2.2]{ACSS21}, possibly replacing $X$ by a foliated resolution,

% we may assume that $\Ff$ is induced by a contraction $f: X\rightarrow Z$. Possibly blowing up more, we may assume that $S$ is normal, and the induced morphism \hl{$f|_S: S\rightarrow f(S)$ is equi-dimensional}\footnote{Om: It may be difficult to keep $f|_S$ equidimensional, $S$ and $X$ normal all these at the same time.}. 

% If $S$ is not $\Ff$-invariant, then $f(S)=Z$, and $\Ff_S$ is the foliation induced by $f|_S$. In particular, (1) holds. Since the leaves of $\Ff_S$ are intersection of $S$ and the leaves of $f$, we get (2).

Now we assume that $S$ is $\Ff$-invariant. From the above computation we see that it is enough to prove the statement for $\Ff_{S'}$ on $X'$. Thus replacing $f:X\to Z$ by $f':X'\to Z'$, $S$ by $S'$, and $\Ff$ by $\Ff'$, we may assume that $\Ff$ is induced by an equidimensional contraction $f:X\to Z$. 
We pick an open subset $U\subset X$ such that $U$ contains the generic point of $S$ and does not intersect $\Sing(\Ff)\cup\Sing(X)\cup\Sing(S)$, and let $S^0:=S\cap U$. Then the leaves of the foliation on $S^0$ induced by the natural inclusion $T_{\Ff}|_{S^0}\rightarrow T_{X}|_{S^0}$ are exactly the intersection of $S^0$ and the the fibers of $f$. Therefore, $\Ff_S$ is induced by the morphism $f|_S\circ \pi: S^\nu\rightarrow f(S)$, so $\Ff_S$ is algebraically integrable, which implies (1). Moreover, the general leaves of $\Ff_S$ are also general leaves of $\Ff$, which implies (2).
\end{proof}

\subsection{Cutting foliations by general hyperplane sections}

One common strategy to deal with usual varieties and pairs is to cut them by general hyperplane sections. 

\begin{defn}
Let $X$ be a normal quasi-projective variety. A divisor $H$ on $X$ is called a \emph{general hyperplane section} if $H$ is very ample and is a general member of the linear system $|H|$.
\end{defn}

The strategy of cutting varieties by general hyperplane sections usually can preserve the singularity of pairs (cf. \cite[Lemma 5.17]{KM98}), but usually cannot preserve the singularities of foliations. We may consider the following example:

\begin{ex}
%\footnote{Om: I added a bit more detail here. Liu: I have made some modifications here. The previous argument ``$\{x\}$ is tangent to $\Ff$ so $E$ is $\tilde\Ff$-invariant" does not sound right for me. I think here we need to use the fact that $H$ is general so $H$ does not contain any dicritical singularity of $\Ff$.}
Let $(X,\Ff,B)$ be an lc foliated triple such that $\dim X=2$, $X$ is klt, and $\rk\Ff=1$. Let $H$ be a general hyperplane on $X$. Then every closed point $x\in H$ is a lc center of $(X, \Ff, H)$; indeed, if $\pi:\widetilde X\to X$ is the blow up of $X$ at $x$ with exceptional divisor $E$ and pullback foliation $\widetilde\Ff$, then $E$ is $\widetilde\Ff$-invariant (as $H$ is general, $\Ff$ is non-dicritical near $x$). Moreover, from the discrepancy computation we see that $a(E, \Ff, X, H)=0=-\epsilon_{\Ff}(E)$. Thus for any $x\in \Supp(H)\cap\Supp(B)$, $a(E,\Ff,X,B+H)<0=-\epsilon_{\Ff}(E)$, so $(X, \Ff, B+H)$ is not lc near $x$.

%Let $H$ be a prime divisor on $X$. Then $(X,B+H)$ is not lc at any point $x\in B\cap H$.
\end{ex}

In this paper, we introduce one way to resolve this issue for algebraically integrable foliations. Instead of cutting a foliation by ``general hyperplane sections", we may cut the foliation by some invariant hyperplane sections, which may not be general, and some non-invariant hyperplanes. 

\subsubsection{Cutting by invariant hyperplanes}

First, we show that we can cut foliations by invariant base-point-free linear systems freely.

\begin{prop}\label{prop: general hyperplane invariant}
Let $(X,\Ff,B)$ be a foliated triple and $W$ a proper sub-variety of $X$. Suppose that $\Ff$ is induced by a morphism $f: X\rightarrow Z$ and $W$ is not tangent to $\Ff$. Let $H_Z\subset Z$ be a general hyperplane section. We define $H:=f^*H_Z$ and $K_{\Ff_H}+B_H:=(K_{\Ff}+B)|_H$, where $\Ff_H$ is the restricted foliation of $\Ff$ on $H$.

Then the following hold:
\begin{enumerate}
    \item $H$ intersects $W$.
    \item If $(X,\Ff,B)$ is (sub-)lc, then $(H,\Ff_H,B_H)$ is (sub-)lc.
     \item For any component $D$ of $\Supp B$ such that $D$ intersects $H$ and any component $C$ of $D\cap H$, $\mult_CB_H=\mult_DB$.
\end{enumerate}
\end{prop}
\begin{proof}
We only need to prove the case when $(X,\Ff,B)$ is sub-lc, and the lc case of (2) follows by (3) and the sub-lc case of (2).

By \cite[Theorem 2.2]{ACSS21}, we have a commutative diagram
\begin{center}$\xymatrix{
X'\ar@{->}[r]^h\ar@{->}_g[d] & X\ar@{->}^f[d] \\
Z'\ar@{->}[r]^{h_Z} & Z.
}$
\end{center}
such that $h$ is a foliated log resolution of $(X,\Supp B)$ and $\Ff':=h^{-1}\Ff$ is induced by $g$. We let $K_{\Ff'}+B':=h^*(K_{\Ff}+B)$ and $H':=h^*H$. Since $H_Z$ is a general hyperplane section, $H$ is general in the base-point-free linear system $|H|$. We let $H_{Z'}:=h_Z^*H_Z$ and $K_{\Ff_{H'}}+B_{H'}:=(K_{\Ff'}+B')|_{H'}$.

(1) Since  $W$ is not tangent to $\Ff$,  $W':=h^{-1}(W)$ is not tangent to $\Ff'$. Thus $\dim g(W')\geq 1$, so $H_{Z'}$ intersects $g(W')$, hence $H_Z$ intersects $h_Z(g(W'))=f(W)$. Thus $H$ intersects $W$.

(2) We let $\tilde B':=B'^{\geq 0}$.  By Lemma \ref{lem: foliated log smooth imply lc}, $(X',\Ff',\tilde B')$ is lc.
%\footnote{Om: In the foliated log resolution of \cite[Theorem 2.2]{ACSS21}, $X'$ is not smooth and neither is $B'$ SNC, so I am not sure how to say $(X', \mathcal F',B')$ is lc. Liu: I think we can directly cite \cite[Lemma 3.1]{ACSS21} (which is our Lemma \ref{lem: foliated log smooth imply lc}). Is it good now?}
Let $K_{\Ff_{H'}}+\tilde B_{H'}:=(K_{\Ff'}+\tilde B')|_{H'}$. By \cite[Proposition 3.2]{ACSS21}, $(H',\Ff_{H'},\tilde B_{H'})$ is lc. Since $\tilde B'\geq B'$, $\tilde B_{H'}\geq B_{H'}$, so $(H',\Ff',B_{H'})$ is sub-lc. Since $H'$ is general,
$$K_{\Ff_{H'}}+B_{H'}=(h|_{H'})^*(K_{\Ff_H}+B_H),$$
so $(H,\Ff_H,B_H)$ is sub-lc.

(3) Since $H$ is general, $(X,B+H)$ is log smooth near the generic point of $C$. Let $B=\sum b_iB_i$ where $B_i$ are the irreducible components of $B$, then
$$B_H=B|_H=\sum b_i(B_i\cap H)$$
near the generic point of $C$. Since $H$ is general, there exists a unique index $i$ such that $B_i\cap H\not=0$ at the genetic point of $C$. Then $B_i\cap H=C$, $B_i=D$, hence $\mult_CB_H=b_i=\mult_DB$.
\end{proof}

\subsubsection{Cutting by non-invariant hyperplanes}

Next we show that, if we only consider the local property of foliations, then we can cut foliation by non-invariant hyperplane sections. We first prove a lemma.

\begin{lem}\label{lem: toroidal cut general hyperplane still lc}
    Let $f: (X,\Sigma)\rightarrow (Z,\Sigma_Z)$ be a toroidal morphism and $z\in Z$ a closed point. Let $\Ff$ be the foliation induced by $f$ and let $B$ be the horizontal$/Z$ part of $\Sigma$. Let $H$ be a general member of a base-point-free linear system on $X$, such that $H$ dominates $Z$.  Then $(X,\Ff,B+H)$ is lc over a neighborhood of $z$.
\end{lem}
\begin{proof}
    The proof is almost similar to the proof of \cite[Lemma 3.1]{ACSS21}. 

    Since $H$ is general, $(X,B+H)$ is lc over a neighborhood of $z$.
   
    We let $x\in X$ be a closed point such that $f(x)=z$. Then locally analytically near $x$, we may assume that $(X,\Sigma)$ is toric. Suppose that $\dim Z=q$, then $\Ff$ is defined by a logarithmic torus invariant $q$-form $$\omega=\sum_{|I|=q}\alpha_I\frac{dx_I}{x_I},$$
    for some $\alpha_I\in\mathbb C$. Let $D$ be the polar locus of $\omega$, then
    $$\mathcal{O}_{X}(K_{[X/\Ff]})=\mathcal{O}_X(-D),$$
    and $D$ is $\Ff$-invariant and torus invariant. In particular, since $H$ and $B$ are horizontal$/Z$, $\Supp D$ and $\Supp(B+H)$ do not have a common component.

    Let $h: X'\rightarrow X$ be a birational morphism, $\Ff':=h^{-1}\Ff$, and $\omega':=h^*\omega$. Then $\omega'$ is a logarithmic $q$-form with polar locus equal to $h^{-1}D+\sum (1-\epsilon_{\Ff}(E_i))E_i$, where $E_i$ are the $h$-exceptional prime divisors. Therefore,
    $$K_{[X'/\Ff']}=h^{-1}_*K_{[X/\Ff]}+\sum (1-\epsilon_{\Ff}(E_i))E_i.$$
Suppose that
$$K_{X'}+h^{-1}_*(D+B+H)=h^*(K_X+D+B+H)+F,$$
then
$$K_{\Ff'}+h^{-1}_*(B+H)=h^*(K_\Ff+B+H)+F+\sum (1-\epsilon_{\Ff}(E_i))E_i.$$
Since $H$ is general, $(X,D+B+H)$ is lc over a neighborhood of $z$, $(X,\Ff,B+H)$ is lc near $x$. Since $x$ can be any closed point over $z$, $(X,\Ff,B+H)$ is lc over $z$. Thus $(X,\Ff,B+H)$ is lc over a neighborhood of $z$.
\end{proof}

\begin{prop}\label{prop: general hyperplane non-invariant}
Let $(X,\Ff,B)$ be a foliated triple and $W$ a proper sub-variety of $X$. Suppose that $\Ff$ is algebraically integrable, $W$ is tangent to $\Ff$, and $\dim W\geq 1$. Let $H\subset X$ be a general hyperplane section. Let $K_{\Ff_H}+B_H:=(K_{\Ff}+B+H)|_H$, where $\Ff_H$ is the restricted foliation of $\Ff$ on $H$.

Then the following hold:
\begin{enumerate}
    \item $H$ intersects $W$.
    \item For any component $D$ of $\Supp B$ such that $D$ intersects $H$ and any component $C$ of $D\cap H$, $\mult_CB_H=\mult_DB$.
    \item  If $(X,\Ff,B)$ is (sub-)lc, then $(H,\Ff_H,B_H)$ is (sub-)lc near $W|_H$.
\end{enumerate}
\end{prop}
\begin{proof}
(1) is obvious.\\
(2) By \cite[Proposition 3.6]{Dru21}, $K_{\Ff_H}=(K_{\Ff}+H)|_H$, so $B_H=B|_H$.  Let $B=\sum b_iB_i$ where $B_i$ are the irreducible components of $B$, then
$$B_H=B|_H=\sum b_i(B_i\cap H).$$
Since $H$ is general, there exists a unique index $i$ such that $B_i\cap H\not=0$ at the generic point of $C$. Then $B_i\cap H=C$, $B_i=D$, hence $\mult_CB_H=b_i=\mult_DB$. This implies (2).

(3) We let $h: X'\rightarrow X$ be a foliated log resolution of $(X,\Ff,B)$, $\Ff':=h^{-1}\Ff$, $K_{\Ff'}+B':=h^*(K_{\Ff}+B)$, $\tilde B':=B'^{\Ff'}$, $H':=h^*H$, and $W':=h^{-1}(W)$. Then $(X',\Ff',\tilde B')$ is foliated log smooth. Therefore, there exists a toroidal contraction $f: (X',\Sigma_{X'})\rightarrow (Z,\Sigma_Z)$ such that $\Ff'$ is induced by $f$, $(Z,\Sigma_Z)$ is log smooth, and $\Supp\tilde B'\subset\Supp\Sigma_{X'}$ (see Definition \ref{defn: foliated log smooth}). We let $z$ be the image of $W'$ on $Z$. Since $(X,\Ff,B)$ is lc, $(X',\Ff',\tilde B')$ is lc. Moreover, all components of $\tilde B'$ are horizontal$/Z$. By Lemma \ref{lem: toroidal cut general hyperplane still lc}, $(X',\Ff',\tilde B'+H')$ is lc over a neighborhood of $z$. In  particular, $(X',\Ff',\tilde B')$ is lc near $W'|_{H'}$.

Let $K_{\Ff_{H'}}+\tilde B_{H'}:=(K_{\Ff'}+\tilde B')|_{H'}$ and $K_{\Ff_{H'}}+B_{H'}:=(K_{\Ff'}+B')|_{H'}$. By \cite[Proposition 3.2]{ACSS21}, $(H',\Ff_{H'},\tilde B_{H'})$ is lc near $W'|_{H'}$. Since $\tilde B'\geq B'$, $\tilde B_{H'}\geq B_{H'}$, so $(H',\Ff_{H'},B_{H'})$ is sub-lc near $W'|_{H'}$. Since $H$ is general, $K_{\Ff_{H'}}+B_{H'}=h|_{H'}^*(K_{\Ff_H}+B_H)$, so $(H,\Ff_H,B_H)$ is sub-lc near $W|_H$. Finally, if $(X,\Ff,B)$ is lc, then by (2), $B_H\geq 0$, so $(H,\Ff_H,B_H)$ is lc near $W|_H$.
\end{proof}

\subsection{Precise adjunction formulas}

Now we are going to discuss the precise adjunction formulas for foliations, that is, adjunction formulas with controlled coefficients.

\begin{thm}[Adjunction to invariant divisors]\label{thm: adj to invariant}
Let $(X,\Ff,B=\sum_{j=1}^mb_jB_j)$ be a $\Qq$-factorial lc f-triple such that $\Ff$ is induced by a contraction $\pi: X\rightarrow Z$ and $B_j$ are the irreducible components of $\Supp B$. Let $S$ be an $\Ff$-invariant prime divisor on $X$ and $S^\nu\rightarrow S$ the normalization of $S$. Then there exists a restricted foliation $\Ff_S$ on $S^\nu$, prime divisors $C_1,\dots,C_n,T_1,\dots,T_l$ on $S^\nu$, positive integers $w_1,\dots,w_n$,  and non-negative integers $\{w_{i,j}\}_{1\leq i\leq n,1\leq j\leq m}$ satisfying the following. For any real numbers $b_1',\dots,b_m'$,
\begin{enumerate}
    \item $$\left(K_{\Ff}+\sum_{j=1}^mb_j'B_j\right)\Bigg|_{S^\nu}=K_{\Ff_S}+\sum_{i=1}^n\frac{w_i-1+\sum_{j=1}^mw_{i,j}b_j'}{w_i}C_i+\sum_{i=1}^lT_i.$$
    \item If $(X,\Ff,\sum_{j=1}^mb_j'B_j)$ is lc, then $$\left(S^\nu,\Ff_S,\frac{w_i-1+\sum_{j=1}^mw_{i,j}b_j'}{w_i}C_i+\sum_{i=1}^lT_i\right)$$ is lc.
\end{enumerate}
\end{thm}
\begin{proof}
(2) follows from \cite[Proposition 3.2]{ACSS21} and (1) so we only need to prove (1). To this end, we may shrink $X$ near the generic point of a codimension $1$ point of $S^\nu$ and then cutting by hyperplane sections as in Propositions \ref{prop: general hyperplane invariant} and \ref{prop: general hyperplane non-invariant} we may assume that $X$ is a surface and $S$ is an $\Ff$-invariant curve. Then (1) follows from \cite[Theorem 3.2(1)]{LMX23b} and \cite[Proposition 3.2]{ACSS21}.
\end{proof}

\begin{thm}[Adjunction to non-invariant divisors]\label{thm: adj to non-invariant}
Let $(X,\Ff,B=S+\sum_{j=1}^mb_jB_j)$ be a $\Qq$-factorial lc f-triple such that $\Ff$ is induced by a contraction $\pi: X\rightarrow Z$ and $S,B_j$ are the irreducible components of $\Supp B$. Let $S^\nu\rightarrow S$ the normalization of $S$. Then there exists a restricted foliation $\Ff_S$ on $S^\nu$, prime divisors $C_1,\dots,C_n,T_1,\dots,T_l$ on $S^\nu$, positive integers $w_1,\dots,w_n$,  and non-negative integers $\{w_{i,j}\}_{1\leq i\leq n,1\leq j\leq m}$ satisfying the following. For any real numbers $b_1',\dots,b_m'$,
\begin{enumerate}
    \item $$\left(K_{\Ff}+S+\sum_{j=1}^mb_j'B_j\right)\Bigg|_{S^\nu}=K_{\Ff_S}+\sum_{i=1}^n\frac{w_i-1+\sum_{j=1}^mw_{i,j}b_j'}{w_i}C_i+\sum_{i=1}^lT_i.$$
    \item If $(X,\Ff,S+\sum_{j=1}^mb_j'B_j)$ is lc, then $$\left(S^\nu,\Ff_S,\frac{w_i-1+\sum_{j=1}^mw_{i,j}b_j'}{w_i}C_i+\sum_{i=1}^lT_i\right)$$ is lc.
\end{enumerate}
\end{thm}
\begin{proof}
(2) follows from \cite[Proposition 3.2]{ACSS21} and (1) so we only need to prove (1). Thus shrinking $X$ near the generic point of a codimension $1$ point of $S^\nu$ and then cutting general members of a by base-point-free linear system as in Propositions \ref{prop: general hyperplane invariant} and \ref{prop: general hyperplane non-invariant}, and we may assume that $X$ is a surface and $S$ is a non-$\Ff$-invariant curve. Thus either $\rk\Ff=2$, or $\rk\Ff=1$. If $\rk\Ff=1$, then by the classification of foliated lc surface singularities (cf.  \cite[Theorem 3.19]{LMX23a}), $S$ does not intersect $B_j$ for any $j$ and $K_{\Ff}+S$ is Cartier near $x$, and the theorem follows. If $\rk\Ff=2$, then (1) follows from the usual adjunction formula (cf. \cite[Theorem 3.10]{HLS19}).
\end{proof}

\section{Property (*) modification and ACSS modification}\label{sec: acss modification}

\cite{ACSS21} famously introduced the concepts of Property $(*)$ pairs \cite[Definition 2.13]{ACSS21}, Property $(*)$ modifications, and Property $(*)$ foliations, and use these new concepts to prove the cone theorem for algebraically integrable foliations \cite[Theorem 3.9]{ACSS21}. In \cite{CS23} the authors have further developed the theory of Property $(*)$ modifications. 

We expect that Property $(*)$ modifications for foliations should be similar to dlt modifications for usual pairs. However, in the formal definition of Property $(*)$ modifications \cite[Definition 3.8]{ACSS21}, some important properties, e.g. $\Qq$-factoriality, are missing. \cite{CS23} resolves this issue by defining a new type of Property $(*)$ modifications. However, when running the MMP for Property $(*)$ foliations, there are a lost of repetitive arguments on some additional properties of Property $(*)$ foliations (e.g. the equidimensional contraction $X\rightarrow Z$ and the associated divisor $G$, cf. \cite[Proposition 2.3]{CS23}). 

To avoid those repetitive argument in this section, we shall introduce a new class of foliated triples, and we call them \emph{ACSS foliated triples}, in honor of Ambro-Cascini-Shokurov-Spicer. Roughly speaking, these triples are foliated triples satisfying Property $(*)$ as well as the nice properties mentioned above. We will show that this class of foliated triples is preserved under the MMP, and show the existence of an ``ACSS modification" for algebraically integrable foliations, i.e. a Property $(*)$ modification that is $\Qq$-factorial and also satisfies the ``ACSS property", as introduced in Definition \ref{defn: ACSS f-triple}.

\subsection{ACSS foliated triples and some basic properties}
\begin{defn}\label{defn: pair property *}
Let $(X,B)/Z$ be a sub-lc sub-pair and the induced morphism $\pi: X\rightarrow Z$ is a contraction. Let $B^v$ and $B^h$ be the vertical$/Z$ and horizontal$/Z$ part of $B$ respectively. We say that $(X,B)/Z$ satisfies \emph{Property $(*)$} if
\begin{enumerate}
    \item there exists a reduced divisor $\Sigma_Z$ on $Z$, such that $(Z,\Sigma_Z)$ is log smooth and $B^v=\pi^{-1}(\Sigma_Z)$, and
    \item for any closed point $z\in Z$ and any reduced divisor $\Sigma\geq\Sigma_Z$ such that $(Z,\Sigma)$ is log smooth near $z$, $(X,B+\pi^*(\Sigma-\Sigma_Z))$ is sub-lc over a neighborhood of $z$.
\end{enumerate}
Note that by \cite[Lemma 2.14(1)]{ACSS21} our definition is equivalent to \cite[Definition 2.13]{ACSS21}.
\end{defn}

\begin{defn}\label{defn: foliation property *}
Let $(X,\Ff,B)$ be an f-sub-triple. 

Let $G\geq 0$ be a reduced divisor on $X$ and let $\pi: X\rightarrow Z$ be a projective morphism. We say that $(X,\Ff,B;G)/Z$ satisfies Property $(*)$ if the following holds:
\begin{enumerate}
\item $(X,B+G)/Z$ satisfies Property $(*)$. In particular, $\pi$ is a contraction.
\item $\Ff$ is induced by $\pi$.
\item $G$ is an $\Ff$-invariant divisor.
\end{enumerate}
We say that $(X,\Ff,B)/Z$ satisfies Property $(*)$ if $(X,\Ff,B;G)/Z$ satisfies Property $(*)$ for some reduced $\Ff$-invariant divisor $G\geq 0$. We say that $(X,\Ff,B)$ satisfies Property $(*)$ if $(X,\Ff,B;G)/Z$ satisfies Property $(*)$ for some $G$ and $X\rightarrow Z$. We say that $\pi: X\rightarrow Z$ is an \emph{associated contraction} of $(X,\Ff,B)$ and $Z$ an \emph{associated base} of $(X,\Ff,B)$. We also say that $\pi: X\rightarrow Z$ and $Z$ are \emph{associated to} $(X,\Ff,B)$. For any divisor $G$ such that $(X,\Ff,B;G)/Z$ satisfies Property $(*)$, we say that $G$ is \emph{associated to} $(X,\Ff,B)/Z$.

We remark that the choice of $\pi$ and $G$ may not be unique.
\end{defn}

\begin{defn}\label{defn: ACSS f-triple}
Let $(X,\Ff,B)$ be an f-triple, $G\geq 0$ a reduced divisor on $X$, and $\pi: X\rightarrow Z$ a projective morphism. We say that $(X,\Ff,B;G)/Z$ is \emph{weak ACSS} if 
\begin{enumerate}
    \item $(X,\Ff,B;G)/Z$ satisfies Property $(*)$ and $(X,\Ff,B)$ is lc, and
    \item $\pi$ is equidimensional.
\end{enumerate}
In this case, we say that the divisor $G$ and the variety $Z$ are \emph{properly associated to} $(X,\Ff,B)$. We say that $(X,\Ff,B;G)/Z$ is \emph{ACSS} if the following additional conditions are satisfied:
\begin{enumerate}
    \item[(3)] for any divisor $\Sigma$ on $Z$ such that $\Sigma\geq\pi(G)$ and $(Z,\Sigma)$ is log smooth, $$(X,B+G+\pi^*(\Sigma-\pi(G)))$$ is qdlt (cf. \cite[Definition 35]{dFKX17}), and
   \item[(4)] for any lc center $W$ of $(X,\Ff,B)$ with generic point $\eta_W$, $(X,B+G)$ is toroidal near $\eta_W$, $\pi$ is a toroidal morphism near $\eta_W$, and $W$ is an lc center of $(X,\Ff,\lfloor B\rfloor)$ (in particular, $K_{\Ff}+\lfloor B\rfloor$ is always $\Qq$-Cartier near $\eta_W$).
  \end{enumerate}
%then we say that $(X,\Ff,B;G)/Z$ is \emph{ACSS}, and we say that $G$ is \emph{properly associated} to $(X,\Ff,B)/Z$.

We say that $(X,\Ff,B)/Z$ is \emph{weak ACSS} (resp. \emph{ACSS}) if $(X,\Ff,B;G)/Z$ is \emph{weak ACSS} (resp. \emph{ACSS}) for some $G$. We say that $(X,\Ff,B)$ is \emph{weak ACSS} (resp. \emph{ACSS}) if $(X,\Ff,B;G)/Z$ is \emph{weak ACSS} (resp. \emph{ACSS}) for some $G$ and $X\rightarrow Z$.

We remark that the condition ``weak ACSS" relies on the choice of $\pi: X\rightarrow Z$ but does not rely on the choice of $G$. That is, for any $(X,\Ff,B;G)/Z$ and $(X,\Ff,B;G')/Z$ which satisfy Property $(*)$, $(X,\Ff,B;G)/Z$ is weak ACSS if and only if $(X,\Ff,B;G')/Z$ is weak ACSS. On the other hand, the condition ``ACSS" relies on both the choice of $\pi: X\rightarrow Z$ and the choice of $G$. That is, it is possible that  $(X,\Ff,B;G)/Z$ and $(X,\Ff,B;G')/Z$ both satisfy Property $(*)$, but $(X,\Ff,B;G)/Z$  is ACSS while $(X,\Ff,B;G')/Z$ is not.
\end{defn}

The following lemmas will be very useful when applying to the minimal model program for algebraically integrable foliations.

\begin{lem}\label{lem: acss smaller coefficient}
Let $(X,\Ff,B)$ be an f-triple satisfying Property $(*)$ with an associated contraction $X\rightarrow Z$. Let $G$ be a divisor associated to $(X,\Ff,B)/Z$. Let $B'$ be an $\Rr$-divisor on $X$ such that $K_{\Ff}+B'$ is $\Rr$-Cartier and $B\geq B'\geq 0$. Then:
\begin{enumerate}
\item $(X,\Ff,B';G)/Z$ satisfies Property $(*)$.
\item If $(X,\Ff,B;G)/Z$ is weak ACSS, then $(X,\Ff,B';G)/Z$ is weak ACSS.
\item If $(X,\Ff,B;G)/Z$ is ACSS, then $(X,\Ff,B';G)/Z$ is ACSS.
\end{enumerate}
\end{lem}
\begin{proof}
(1) We check the conditions of Definition \ref{defn: foliation property *} for $(X,\Ff,B';G)/Z$. To check Definition \ref{defn: foliation property *}(1), we need to check Definition \ref{defn: pair property *} for $(X,B'+G)/Z$. Since $B\geq 0$ and $(X,\Ff,B)$ is lc, all components of $B$ are horizontal$/Z$, so the vertical$/Z$ part of $B+G$ and $B'+G$ are equal. This verifies Definition \ref{defn: pair property *}(1) for $(X,B'+G)/Z$. Definition \ref{defn: pair property *}(2) holds for $(X,B'+G)/Z$ because $B\geq B'$.  Definition \ref{defn: foliation property *}(2)(3) are clear since they do not rely on $B'$. 

(2) We check (1)(2) of Definition \ref{defn: ACSS f-triple} for $(X,\Ff,B';G)/Z$. By (1), $(X,\Ff,B';G)/Z$ satisfies Property (*). Since $B\geq B'$, $(X,\Ff,B')$ is lc. Thus Definition \ref{defn: ACSS f-triple} is verified for $(X,\Ff,B';G)/Z$. Definition \ref{defn: ACSS f-triple}(2) does not rely on $B'$, so it is clear.

(3) We check (3)(4) of Definition \ref{defn: ACSS f-triple} for $(X,\Ff,B';G)/Z$. Since $B\geq B'$, (3) of  Definition \ref{defn: ACSS f-triple} for $(X,\Ff,B';G)/Z$ is verified. For any lc center $W$ of $(X,\Ff,B')$ with generic point $\eta_W$, $W$ is an lc center of $(X,\Ff,B)$, hence $(X,B+G)$ is toroidal near $\eta_W$, $\pi$ is a toroidal morphism near $\eta_W$, and $W$ is an lc center of $(X,\Ff,\lfloor B\rfloor)$. Since $B\geq B'$, $(X,B'+G)$ is toroidal near $\eta_W$. Moreover, since $(X,\Ff,B)$ and $(X,\Ff,B')$ are both lc, $\eta_W$ is not contained in $\Supp(B-B')$. Since
$$\Supp(\lfloor B\rfloor-\lfloor B'\rfloor)=\Supp(\lfloor B'+(B-B')\rfloor-\lfloor B'\rfloor)\subset\Supp (B-B'),$$
$\eta_W$ is not contained in $\Supp(\lfloor B\rfloor-\lfloor B'\rfloor)$. Thus $W$ is an lc center of $(X,\Ff,\lfloor B'\rfloor)$. This verifies (4) of  Definition \ref{defn: ACSS f-triple} and we are done.
\end{proof}

\begin{lem}\label{lem: alg int foliation lct achieved}
Let $(X,\Ff,B)$ be a sub-lc f-sub-triple and $D$ an $\Rr$-Cartier $\Rr$-divisor on $X$. Suppose that $\Ff$ is algebraically integrable. Let
$$t:=\sup\{s\mid s\geq 0, (X,\Ff,B+sD)\text{ is sub-lc}.\}$$
Then $t=+\infty$ or $t=\max\{s\mid s\geq 0, (X,\Ff,B+sD)\text{ is sub-lc}\}.$ In particular, $(X,\Ff,B+tD)$ is \text{sub-lc} if $t<+\infty$.
\end{lem}
\begin{proof}
By \cite[Theorem 2.2]{ACSS21}, we may assume that $(X,\Ff,\Supp B\cup\Supp D))$ is foliated log smooth. By Lemma \ref{lem: foliated log smooth imply lc},
\begin{align*}
    t&=\sup\{s\mid s\geq 0, a(E,X,\Ff,B+sD)\geq-\epsilon_{\Ff}(E)\text{ for any prime divisor } E\text{ on }X\}\\
    &=\sup\{s\mid s\geq 0, a(E,X,\Ff,B+sD)\geq-\epsilon_{\Ff}(E)\text{ for any prime divisor } E\subset\Supp D\}.
\end{align*}
Since there are only finitely many components of $\Supp D$, the lemma follows.
\end{proof}

\begin{lem}\label{lem: acss f-triple perturb coefficient}
Let $(X,\Ff,B)$ be a f-triple, $G\geq 0$ a reduced divisor on $X$, and $\pi: X\rightarrow Z$ a contraction, such that $(X,\Ff,B;G)/Z$ is ACSS. Let $D\geq 0$ be an $\Rr$-Cartier $\Rr$-divisor on $X$ such that $\Supp D\subset\Supp\{B\}$. Then there exists a positive real number $\delta$, such that
$(X,\Ff,B+\delta D;G)/Z$ is ACSS. In particular, $(X,\Ff,B+\delta D)/Z$ and $(X,\Ff,B+\delta D)$ are ACSS.
\end{lem}
\begin{proof}
We check the conditions (1-4) of Definition \ref{defn: ACSS f-triple} for $(X,\Ff,B+\delta D;G)/Z$. (2) of Definition \ref{defn: ACSS f-triple} is clear. For any divisor $\Sigma$ on $X$ such that $\Sigma\geq\pi(G)$ and $(Z,\Sigma)$ is log smooth, $$(X,B+G+\pi^*(\Sigma-\pi(G)))$$ is qdlt. Thus, there exists a positive real number $\delta$ such that $(X,B+\delta D+G+\pi^*(\Sigma-\pi(G)))$ is qdlt. This implies (3) of Definition \ref{defn: ACSS f-triple}, and shows that $(X,B+\delta D;G)/Z$ satisfies Property $(*)$. 

\begin{claim}\label{claim: D does not contain lc center}
$\Supp D$ does not contain any lc center of $(X,\Ff,B)$.
\end{claim}
\begin{proof}
Suppose not. Then there exists a prime divisor $E$ over $X$ such that $a(E,\Ff,B)=-\epsilon_{\Ff}(E)$ and $\mult_ED>0$. Let $W:=\Center_XE$ and let $\eta_W$ be the generic point of $W$. Since $(X,\Ff,B;G)/Z$ is ACSS, $(X,B+G)$ is toroidal near $\eta_W$, $\pi$ is a toroidal morphism near $\eta_W$, and $W$ is an lc center of $(X,\Ff,\lfloor B\rfloor)$. Thus $\eta_W$ is not contained in $\Supp\{B\}$, so $\eta_W$ is not contained in $\Supp D$, a contradiction. 
\end{proof}

\noindent\textit{Proof of Lemma \ref{lem: acss f-triple perturb coefficient} continued}. By Claim \ref{claim: D does not contain lc center} and  Lemma \ref{lem: alg int foliation lct achieved}, there exists a positive real number $\delta$ such that $(X,\Ff,B+\delta D)$ is lc. This implies (1) of Definition \ref{defn: ACSS f-triple}. 

Possibly replacing $\delta$ with a smaller positive real number, we may assume that any lc center of $(X,\Ff,B+\delta D)$ is an lc center of $(X,\Ff,B)$. In particular, $\lfloor B+\delta D\rfloor=\lfloor B\rfloor$. By Claim \ref{claim: D does not contain lc center}, $\Supp D$ does not contain any lc center of $(X,\Ff,B+\delta D)$. Since $(X,\Ff,B;G)/Z$ is ACSS, for any lc center $W$ of $(X,\Ff,B)$,
\begin{itemize}
    \item $(X,B+G)$ is toroidal near $\eta_W$, and since $\Supp D\subset\Supp\{B\}$, $(X,B+\delta D+G)$ is toroidal near $\eta_W$,
    \item $\pi$ is a toroidal morphism near $\eta_W$, and
    \item $W$ is an lc center of $(X,\Ff,\lfloor B\rfloor)=(X,\Ff,\lfloor B+\delta D\rfloor)$.
\end{itemize}
Thus (4) of Definition \ref{defn: ACSS f-triple} holds for  $(X,\Ff,B+\delta D;G)/Z$, and we are done.
\end{proof}

\subsection{Minimal model program for ACSS foliated triples}

\begin{lem}\label{lem: ACSS mmp can run}
Let $(X,\Ff,B)$ be an f-triple, $G\geq 0$ a reduced divisor on $X$, and $\pi: X\rightarrow Z$ a contraction, such that $(X,\Ff,B;G)/Z$ is weak ACSS. Then:
\begin{enumerate}
    \item Any $(K_{\Ff}+B)$-negative extremal ray is a $(K_X+B+G)$-negative extremal ray$/Z$.
    \item We may run a step of a $(K_{\Ff}+B)$-MMP, and any such step is a step of a $(K_X+B+G)$-MMP$/Z$.
    \item For any step $\phi: (X,\Ff,B)\dashrightarrow (X',\Ff',B')$ of a $(K_{\Ff}+B)$-MMP that is not a Mori fiber space,
    \begin{enumerate}
    \item $(X',\Ff',B';\phi_*G)/Z$ is weak ACSS,
        \item if $(X,\Ff,B;G)/Z$ is ACSS, then $(X',\Ff',B';\phi_*G)/Z$ is ACSS, and
        \item if $X$ is $\Qq$-factorial, then $X'$ is $\Qq$-factorial.
    \end{enumerate}
    \item Any sequence of steps of a $(K_{\Ff}+B)$-MMP is a sequence of steps of a $(K_X+B+G)$-MMP$/Z$.
\end{enumerate}
\end{lem}
\begin{proof}
(1) follows from \cite[Theorem 3.9]{ACSS21}. (2) follows from (1) and the cone theorem, contraction theorem (cf. \cite{Amb03,Fuj17}), and the existence of flips for lc pairs \cite{Bir12,HX13}. (3.a) follows from (2) and \cite[Proposition 2.18]{ACSS21}. 

If $(X,\Ff,B;G)/Z$ is ACSS, then for any divisor $\Sigma$ on $Z$ such that $\Sigma\geq\pi(G)$ and $(Z,\Sigma)$ is log smooth, $(X,B+G+\pi^*(\Sigma-\pi(G)))$ is qdlt. By (2), $\phi$ is a step of a $(K_X+B+G)$-MMP$/Z$, hence a step of a  $(K_X+B+G+\pi^*(\Sigma-\pi(G)))$-MMP$/Z$. Thus $(X',B'+\phi_*G+\pi'^*(\Sigma-\pi(G)))$ is qdlt (cf. \cite[Paragraph after Definition 35]{dFKX17}), where $\pi': X'\rightarrow Z$ is the induced contraction. Moreover, for any lc place $E$ of $(X',\Ff',B')$, since $$-\epsilon_{\Ff}(E)\leq a(E,\Ff,B)\leq a(E,\Ff',B')\leq -\epsilon_{\Ff'}(E)=-\epsilon_{\Ff}(E),$$
$E$ is also an lc place of $(X,\Ff,B)$ and $\phi$ is an isomorphism near the generic point of $\Center_XE$. Thus (3.b) follows from (3.a). (3.c) follows from (3.a) and \cite[Corollary 5.20, Theorem 6.2]{HL21} (which are essentially the same lines of the proofs of \cite[Corollaries 3.17, 3.18]{KM98}).

(4) follows from (2) and (3).
\end{proof}

\begin{nota}
Let $(X_0,\Ff_0,B_0)$ be a weak ACSS f-triple, $Z$ a base properly associated to $(X_0,\Ff_0,B_0)$, and $G_0$ a divisor associated to $(X_0,\Ff_0,B_0)/Z$. When we say the following
\begin{center}$\xymatrix{
(X_0,\Ff_0,B_0;G_0)\ar@{-->}[r]^{f_0} & (X_1,\Ff_1,B_1;G_1)\ar@{-->}[r]^{\ \ \ \ \ \ \ \ f_1} & \dots\ar@{-->}[r] & (X_n,\Ff_n,B_n;G_n)\ar@{-->}[r]^{\ \ \ \ \ \ \ \ f_n} & \dots 
}$
\end{center}
is a (possibly infinite) sequence of steps of a $(K_{\Ff_0}+B_0)$-MMP, we mean the following: for any $i$, $f_i: X_{i}\rightarrow X_{i+1}$ is a step of a $(K_{\Ff_i}+B_i)$-MMP that is not a Mori fiber space, $\Ff_{i+1}:=(f_i)_*\Ff_i$, $B_{i+1}:=(f_i)_*B_i$, and $G_{i+1}:=(f_i)_*G_i$. By Lemma \ref{lem: ACSS mmp can run},
\begin{enumerate}
    \item $(X_i,\Ff_i,B_i;G_i)/Z$ is weak ACSS,
    \item if $(X_0,\Ff_0,B_0;G_0)/Z$ is ACSS, then $(X_i,\Ff_i,B_i;G_i)/Z$ is ACSS, and
    \item if $X_0$ is $\Qq$-factorial, then $X_i$ is $\Qq$-factorial.
\end{enumerate}
\end{nota}

\begin{defn}[Models]\label{defn: models}
Let $(X,\Ff,B)/U$ be an lc f-triple, $\phi: X\dashrightarrow X'$ a birational map over $U$, $E:=\Exc(\phi^{-1})$ the reduced $\phi^{-1}$-exceptional divisor, $\Ff':=\phi_*\Ff$, and $B':=\phi_*B+E^{\Ff'}$.
\begin{enumerate}
    \item $(X',\Ff',B')/U$ is called a \emph{log birational model} of $(X,\Ff,B)/U$. 
    \item $(X',\Ff',B')/U$ is called a \emph{weak lc model} of $(X,\Ff,B)/U$ if 
\begin{enumerate}
\item $(X',\Ff',B')/U$ is a log birational model of $(X,\Ff,B)/U$, 
    \item $K_{\Ff'}+B'$ is nef$/U$, and
    \item for any prime divisor $D$ on $X$ which is exceptional over $X'$, $a(D,\Ff,B)\leq a(D,\Ff',B')$.
\end{enumerate}
\item $(X',\Ff',B')/U$ is called a \emph{log minimal model} of $(X,\Ff,B)/U$ if
\begin{enumerate}
    \item $(X',\Ff',B')/U$ is a weak lc model of $(X,\Ff,B)/U$,
    \item $(X',\Ff',B')$ is $\Qq$-factorial ACSS.
    \item for any prime divisor $D$ on $X$ which is exceptional over $X'$, $a(D,\Ff,B)<a(D,\Ff',B')$.
\end{enumerate}
\item $(X',\Ff',B')/U$ is called a \emph{good minimal model} of $(X,\Ff,B)/U$ if
\begin{enumerate}
        \item $(X',\Ff',B)/U$ is a log minimal model of $(X,\Ff,B)/U$, and
        \item $K_{\Ff'}+B'$ is semi-ample$/U$.
\end{enumerate}
\end{enumerate}
\end{defn}

\begin{thm}\label{thm: mmp very exceptional alg int fol}
Let $(X,\Ff,B)$ be a $\Qq$-factorial weak ACSS f-triple, $h: X\rightarrow Y$ a contraction, and $E\geq 0$ an $\Rr$-divisor on $X$, such that $K_{\Ff}+B\sim_{\mathbb R,Y}E$ and $E$ is very exceptional over $Y$. Then:
\begin{enumerate}
    \item We may run a $(K_{\Ff}+B)$-MMP$/Y$ with scaling of an ample$/Y$ $\Rr$-divisor.
    \item Any $(K_{\Ff}+B)$-MMP$/Y$ with scaling of an ample$/Y$ $\Rr$-divisor terminates with a weak lc model $(X',\Ff',B')/Y$ of $(X,\Ff,B)/Y$ satisfying the following:
    \begin{enumerate}
    \item $K_{\Ff'}+B'\sim_{\mathbb R,Y}0$.
        \item The divisors contracted by the induced birational map $X\dashrightarrow X'$ are exactly $\Supp E$.
        \item If $(X,\Ff,B)$ is ACSS, then $(X',\Ff',B')/Y$ is a good minimal model of $(X,\Ff,B)$. 
    \end{enumerate}
\end{enumerate}
\end{thm}
\begin{proof}
(1) It follows from Lemma \ref{lem: ACSS mmp can run}(1) and the MMP with scaling for lc pairs.

(2) Let $Z$ be a properly associated base of $(X,\Ff,B)$ and let $G$ be a divisor associated to $(X,\Ff,B)/Z$, such that if $(X,\Ff,B)$ is ACSS, then $G$ is properly associated to $(X,\Ff,B)/Z$. By (1), we may run a $(K_{\Ff}+B)$-MMP$/Y$ with scaling of an ample$/Y$ $\Rr$-divisor $A$. This MMP is also a $(K_{\Ff}+B)$-MMP$/Y$ with scaling of $A+h^*H$ where $H$ is a sufficiently ample divisor on $Y$. Thus, possibly replacing $A$ with  $A+h^*H$, we may assume that $A$ is ample. By Lemma \ref{lem: ACSS mmp can run}(1), this  $(K_{\Ff}+B)$-MMP$/Y$ with scaling of $A$ is also a $(K_X+B+G)$-MMP$/Z$ with scaling of $A$, and consists of a sequence of divisorial contractions and flips
$$(X,\Ff,B;G):=(X_1,\Ff_1,B_1;G_1)\dashrightarrow (X_2,\Ff_2,B_2;G_2)\dashrightarrow\dots\dashrightarrow (X_n,\Ff_n,B_n;G_n)\dashrightarrow\dots,$$
such that for any $n$, $(X_n,\Ff_n,B_n;G_n)/Z$ is $\Qq$-factorial weak ACSS, and if $(X,\Ff,B;G)/Z$ is ACSS, then $(X_n,\Ff_n,B_n;G_n)/Z$ is ACSS. For each $i$, we let $A_i,E_i$ be the images of $A,E$ on $X_i$ respectively, and let $\lambda_i$ be the $i$-th scaling number, i.e. $$\lambda_i:=\inf\{t\geq 0\mid K_{\Ff_i}+B_i+tA_i\text{ is nef$/Y$}\}.$$ 

\begin{claim}\label{claim: very exceptional mmp terminates}
The MMP $X_1\dashrightarrow X_2\dashrightarrow\dots\dashrightarrow X_n\dashrightarrow\dots$ terminates. 
\end{claim}
\begin{proof}
Suppose that the MMP does not terminate.  Let $\lambda:=\lim_{i\rightarrow+\infty}\lambda_i$. Since $X_1\dashrightarrow X_2\dashrightarrow\dots\dashrightarrow X_n\dashrightarrow\dots$ is also a $(K_{X}+B+G)$-MMP$/Z$ with scaling of $A$, if $\lambda>0$, then this MMP is also a $(K_{X}+B+G+\lambda A)$-MMP$/Z$, which terminates by \cite[Theorem 1.5]{HH20} and \cite[Theorem 1.9]{Bir12}, a contradiction. Thus $\lambda=0$.

We let $n$ be a positive integer such that the induced birational map $X_i\dashrightarrow X_{i+1}$ is a flip for any $i\geq n$. Then $K_{\Ff_n}+B_n$ is the limit of movable$/Y$ $\Rr$-Cartier $\Rr$-divisors $K_{\Ff_n}+B_n+\lambda_iA_i$.

Since $E$ is very exceptional$/Y$, $E_n$ is very exceptional$/Y$. Therefore, for any prime divisor $S$ on $X_n$, $(K_{\Ff_n}+B_n)\cdot C=E_n\cdot C\geq 0$ for very general curves $C$ of $S/Y$. Since $E_n\geq 0$, by \cite[Lemma 3.3]{Bir12}, $E_n=0$. Thus the MMP terminates, a contradiction.
\end{proof}

\noindent\textit{Proof of Theorem \ref{thm: mmp very exceptional alg int fol} continued}. By Claim \ref{claim: very exceptional mmp terminates}, the MMP terminates at $X_n$ for some positive integer $n$. Then $E_n\geq 0$ is very exceptional$/Y$. By \cite[Lemma 3.3]{Bir12}, $E_n=0$. We have $K_{\Ff_n}+B_n\sim_{\mathbb R,Y}E_n=0$. Let $(X',\Ff',B')=(X_n,\Ff_n,B_n)$, then we get (2.a). In particular, $\Supp E$ are contracted by $X\dashrightarrow X_n$. Since $X\dashrightarrow X_n$ is also an $E$-MMP, any divisor which is not contained in $\Supp E$ is not contracted by $X\dashrightarrow X_n$, which implies (2.b). 

Since the induced birational map $X\dashrightarrow X_n$ does not extract any divisor, $(X_n,\Ff_n,B_n)/Y$ is a log birational model of $(X,\Ff,B)/Y$. Since  $K_{\Ff_n}+B_n\sim_{\mathbb R,Y}0$, $K_{\Ff_n}+B_n$ is semi-ample$/Y$. Since $X\dashrightarrow X_n$ is a sequence of steps of a $(K_{\Ff}+B)$-MMP, $a(D,\Ff,B)<a(D,\Ff_n,B_n)$ for any prime divisor $D$ on $X$ which is exceptional over $X_n$. Thus $(X_n,\Ff_n,B_n)/Y$ is a weak lc model of $(X,\Ff,B)/Y$. Moreover, if $(X,\Ff,B)$ is ACSS, then by Lemma \ref{lem: ACSS mmp can run}(3.b), $(X_n,\Ff_n,B_n)$ is ACSS, hence $(X_n,\Ff_n,B_n)/Y$ is a good minimal model of $(X,\Ff,B)/Y$. This implies (2.c) and we are done.
\end{proof}

\subsection{ACSS modifications}

\begin{defn}
Let $(X,\Ff,B)$ be an lc f-triple such that $\Ff$ is algebraically integrable.
\begin{enumerate}
    \item A \emph{Property $(*)$ modification} of $(X,\Ff,B)$ is a projective birational morphism $f: Y\rightarrow X$, such that
\begin{enumerate}
    \item $$\left(Y,\Ff_Y:=f^{-1}\Ff,B_Y:=f^{-1}_*B+(\Supp\Exc(f))^{\Ff_Y}\right)$$
is a weak ACSS f-triple, 
\item $Y$ is klt, and
\item  $K_{\Ff_Y}+B_Y=f^*(K_{\Ff}+B)$.
\end{enumerate}
In this case, we say that $(Y,\Ff_Y,B_Y)$ is a \emph{Property $(*)$ model} of $(X,\Ff,B)$. If $\pi: Y\rightarrow Z$ is a properly associated contraction of $(Y,\Ff_Y,B_Y)$, then we also say that $(Y,\Ff_Y,B_Y)/Z$ is a \emph{Property $(*)$ model} of $(X,\Ff,B)$. In addition, if
\begin{enumerate}
    \item[(d)] $Y$ is $\Qq$-factorial, 
\end{enumerate}
then we say that $f$ is a \emph{$\Qq$-factorial Property $(*)$ modification} of $(X,\Ff,B)$, and say that $(Y,\Ff_Y,B_Y)/Z$ is a \emph{$\Qq$-factorial Property $(*)$ model} of $(X,\Ff,B)$. 
\item An \emph{ACSS modification} of $(X,\Ff,B)$ is a Property $(*)$ modification $f: Y\rightarrow X$ satisfying the following: let $\Ff_Y:=f^{-1}\Ff$ and $K_{\Ff_Y}+B_Y:=f^*(K_{\Ff}+B)$, then $(Y,\Ff_Y,B_Y)$ is $\Qq$-factorial ACSS.

In this case, we say that $(Y,\Ff_Y,B_Y)$ is an \emph{ACSS model} of $(X,\Ff,B)$.  If $\pi: Y\rightarrow Z$ is a properly associated contraction of $(Y,\Ff_Y,B_Y)$ and $G$ is a divisor properly associated to $(Y,\Ff_Y,B_Y)/Z$, then we also say that $(Y,\Ff_Y,B_Y)/Z$ and $(Y,\Ff_Y,B_Y;G)/Z$ are \emph{ACSS models} of $(X,\Ff,B)$.
\item Let $f: Y\rightarrow X$ be a $\Qq$-factorial Property $(*)$ modification (resp. ACSS modification) of $(X,\Ff,B)$. $\Ff_Y:=f^{-1}\Ff$, and $K_{\Ff_Y}+B_Y:=f^*(K_{\Ff}+B)$. We say that $f$ is a \emph{proper Property $(*)$ modification} (resp. \emph{proper ACSS modification}) of $(X,\Ff,B)$ if there exists a properly associated contraction $Y\rightarrow Z$ of $(Y,\Ff_Y,B_Y)$ and a divisor $G$ associated to (resp. properly associated to) $(Y,\Ff_Y,B_Y)/Z$, such that for any $f$-exceptional $\Ff_Y$-invariant divisor $D$, $D\subset\Supp G$. We call $(Y,\Ff_Y,B_Y)$, $(Y,\Ff_Y,B_Y)/Z$, and $(Y,\Ff_Y,B_Y;G)/Z$ \emph{proper Property $(*)$ models} (resp. \emph{proper ACSS models}) of $(X,\Ff,B)$.
\end{enumerate}
\end{defn}

We remark that our Property $(*)$ modification is the same as \cite[Definition 3.8]{ACSS21}, while the ``Property $(*)$ modification" in \cite[Theorem 2.4]{CS23} corresponds to our proper Property $(*)$ modification.

The following proposition is an analogue of \cite[Theorem 3.10]{ACSS21} and \cite[Theorem 2.4]{CS23}. We usually do not need the full power of the following proposition.

\begin{prop}[Existence of ACSS models]\label{prop: ACSS modification include some divisor}
Let $(X,\Ff,B)$ be a projective lc f-triple such that $\Ff$ is algebraically integrable. Let $E_1,\dots,E_n$ be a prime divisor over $X$ such that $a(E_i,\Ff,B)=-\epsilon_{\Ff}(E_i)$ for each $i$. Then there exists a proper ACSS modification $f: Y\rightarrow X$, such that  $E_1,\dots,E_n$ are on $Y$, i.e. $\Center_YE_i$ is a divisor for each $i$.
\end{prop}
\begin{proof}
The proof is very similar to the proof of \cite[Theorem 2.4]{CS23}. By \cite[Theorem 2.2]{ACSS21}, there exists a foliated log resolution $h: W\rightarrow X$ of $(X,\Ff,B)$. Possibly blowing-up more and applying \cite[Theorem 2.2]{ACSS21} again, we may assume that $E_1,\dots,E_n$ are on $W$. We let $\Ff_W:=h^{-1}\Ff$, $\tilde B_W:=h^{-1}_*B+\Supp\Exc(h)$, $B_W:=h^{-1}_*B+(\Supp \Exc(h))^{\Ff}$, and $\pi: W\rightarrow Z$ an equidimensional toroidal contraction which induces $\Ff_W$.

Since no component of $B_W$ is $\Ff_W$-invariant, by \cite[Proposition 2.16]{ACSS21}, there exists a reduced divisor $G_W\geq 0$ on $W$ such that $G_W$ is vertical$/Z$ and $(W,\Supp B_W+G_W)/Z$ satisfies Property $(*)$. Since $(W,\Ff_W,\tilde B_W)$ is foliated log smooth and any component of $\Exc(h)$ is a component of $\tilde B_W$, possibly adding components to $G_W$, we may assume that all $\Ff_W$-invariant divisors in $\Exc(h)$ are contained in $\Supp G_W$. Since $f$ is a foliated log resolution of $(X,\Ff,B)$, $(W,\Ff_W,\Supp B_W;G_W)/Z$ is $\Qq$-factorial ACSS.

Since
$$K_{\Ff_W}+B_W\sim_{\mathbb R,X}\sum_{F\mid F\text{ is exceptional over }X}(\epsilon_{\Ff}(F)+a(F,\Ff,B))F,$$
by Theorem \ref{thm: mmp very exceptional alg int fol}, we may run a $(K_{\Ff_W}+B_W)$-MMP$/X$ which terminates with a good minimal model $(Y,\Ff_Y,B_Y)/X$ of $(W,\Ff_W,B_W)/X$, such that the divisors contacted by the induced birational map $W\dashrightarrow Y$ are exactly the $h$-exceptional divisors that are not lc places of $(X,\Ff,B)$. In particular, $E_1,\dots,E_n$ are not contracted by the MMP. We may let $f: Y\rightarrow X$ be the induced projective birational morphism and let $B_Y$ be the image of $B_W$ on $Y$. By Lemma \ref{lem: ACSS mmp can run}, $f$ is a proper ACSS modification of $(X,\Ff,B)$. The proposition follows.
\end{proof}

\section{Proof of the main theorems}\label{sec: proof of the main theorems}

In this section, we prove Theorems \ref{thm: acc lct alg int foliation}, \ref{thm: global acc alg int foliation}, \ref{thm: log abundance algebraically integrable foliation}, and Corollary \ref{cor: global acc rank 1 foliation}.

\begin{lem}\label{lem: find nontrivial divisor on ACSS model}
Let $(X,\Ff,B)$ be a projective lc f-triple and $M$ and $\Rr$-Cartier $\Rr$-divisor on $X$, such that 
\begin{itemize}
    \item $\Ff$ is algebraically integrable,
    \item $(X,\Ff,B+M)$ is lc,
    \item $(X,\Ff,B+(1+\epsilon)M)$ is not lc for any positive real number $\epsilon$,
    \item $\Supp B=\Supp(B+M)$, and
    \item for any prime divisor $D$ on $X$ such that $a(D,\Ff,B+M)=-\epsilon_{\Ff}(D)$, $\mult_DM=0$.
\end{itemize}
Then there are two projective birational morphisms $h: X'\rightarrow X$ and $g: Y'\rightarrow X'$ and a real number $t\in (0,1)$ satisfying the following. 
\begin{enumerate}
    \item $h$ is an ACSS modification of $(X,\Ff,B+tM)$ for some $t\in (0,1)$.
    \item For any prime $h$-exceptional divisor $D$, $a(D,\Ff,B)=-\epsilon_{\Ff}(D)$. In particular, $\mult_DM=0$ and $a(D,\Ff,B+sM)=-\epsilon_{\Ff}(D)$ for any real number $s$.
    \item $g$ extracts a unique prime divisor $E$. In particular, $-E$ is ample over $X'$.
    \item $a(E,\Ff,B+M)=-\epsilon_{\Ff}(E)$ and $a(E,\Ff,B)>-\epsilon_{\Ff}(E)$. In particular, $\mult_EM>0$ and $a(E,\Ff,B+sM)>-\epsilon_{\Ff}(E)$ for any real number $s<1$.
    \item Let $B_{Y'},M_{Y'}$ be the strict transforms of $B,M$ on $Y'$ respectively, $\Ff_{Y'}:=(h\circ g)^{-1}\Ff$, and
    $$F_{Y'}:=\sum_{D\mid D\text{ is a prime }h\circ g\text{-exceptional divisor}}\epsilon_{\Ff}(D)D.$$
    Then $(Y',\Ff_{Y'},B_{Y'}+tM_{Y'}+F_{Y'})$ is $\Qq$-factorial ACSS.
\end{enumerate}
\begin{center}$\xymatrix{
Y\ar@{-->}[r]\ar@{->}_f[d] & Y'\ar@{->}^g[d] \\
X & X'\ar@{->}[l]^{h}.
}$
\end{center}
\end{lem}
\begin{proof}
By Lemma \ref{lem: alg int foliation lct achieved}, there exists a prime divisor $P$ on $X$ such that $a(P,\Ff,B+M)=-\epsilon_{\Ff}(P)$ and $\mult_EM\not=0$. By our conditions, $P$ is exceptional over $X$. Since $(X,\Ff,B+M)$ is lc and $\Ff$ is algebraically integrable, by Proposition \ref{prop: ACSS modification include some divisor}, there exists a proper ACSS modification $f: Y\rightarrow X$ of $(X,\Ff,B+M)$ satisfying the following. Let $\Ff_Y:=f^{-1}\Ff$, $B_Y:=f^{-1}_*B,M_Y:=f^{-1}_*M$, and $E_1,\dots,E_n$ the prime $f$-exceptional divisors, then possibly reordering indices, we have $P=E_1$. We let $\pi_Y: Y\rightarrow Z$ be a contraction properly associated to $(Y,\Ff_Y,B_Y+M_Y+\sum_{i=1}^n\epsilon_{\Ff}(E_i)E_i)$, and $G_Y$ a divisor properly associated to $(Y,\Ff_Y,B_Y+M_Y+\sum_{i=1}^n\epsilon_{\Ff}(E_i)E_i)/Z$.

Since $\Supp B=\Supp(B+M)$, by Lemmas \ref{lem: alg int foliation lct achieved} and \ref{lem: acss f-triple perturb coefficient},  $(Y,\Ff_Y,B_Y+tM_Y+\sum_{i=1}^n\epsilon_{\Ff}(E_i)E_i)$ is ACSS for some $t\in (0,1)$. Since $(X,\Ff,B)$ and $(X,\Ff,B+M)$ are lc, $(X,\Ff,B+tM)$ is lc. Thus
$$\left(K_{\Ff_Y}+B_Y+tM_Y+\sum_{i=1}^n\epsilon_{\Ff}(E_i)E_i\right)\sim_{\mathbb R,X}\sum_{i=1}^n\left(\epsilon_{\Ff}(E_i)-a(E_i,\Ff,B+tM)\right)E_i\geq 0.$$
By Theorem \ref{thm: mmp very exceptional alg int fol}, we may run a $(K_{\Ff_Y}+B_Y+tM_Y+\sum_{i=1}^n\epsilon_{\Ff}(E_i)E_i)$-MMP$/X$ which terminates with a good minimal model $(X',\Ff',B'+tM'+F')/X$ of $(Y,\Ff_Y,B_Y+tM_Y+\sum_{i=1}^n\epsilon_{\Ff}(E_i)E_i)/X$ such that $K_{\Ff'}+B'+tM'+F'\sim_{\mathbb R,X}0$, where $B',M',F'$ are the images of $B_Y,M_Y,\sum_{i=1}^n\epsilon_{\Ff}(E_i)E_i$ on $X'$ respectively. Let $G'$ be the image of $G_Y$ on $X'$. By Lemma \ref{lem: ACSS mmp can run}, $(X',\Ff',B'+tM'+F';G')/Z$ is $\Qq$-factorial ACSS. In particular, the induced morphism $h: X'\rightarrow X$ is an ACSS modification of $(X,\Ff,B+tM)$.

By construction, the divisors contracted by the induced birational map $Y\dashrightarrow X'$ are all divisors $E_i$ such that $\epsilon_{\Ff}(E_i)>a(E_i,\Ff,B_Y+tM_Y)$. Since $\mult_{E_1}M\not=0$, $Y\dashrightarrow X'$ contracts $E_1$, hence $Y\dashrightarrow X'$ contains a divisorial contraction. We let $g: Y'\rightarrow X'$ be the last step of the $(K_{\Ff_Y}+B_Y+tM_Y+\sum_{i=1}^n\epsilon_{\Ff}(E_i)E_i)$-MMP$/X$. Since $X'$ is $\Qq$-factorial and $K_{\Ff'}+B'+tM'+F'\sim_{\mathbb R,X}0$, $g$ is a divisorial contraction of a prime divisor $E$. 

We show that $h,g$ and $t$ satisfy our requirements. (1)(5) immediately follow from our construction. For any prime divisor $D$ on $X'$ that is exceptional over $X$,
\begin{align*}
   a(D,\Ff,B+M)&=-\epsilon_{\Ff}(D) && \text{($D$ is also on $Y$ and is exceptional$/X$)}\\
                    &=a(D,\Ff,B+tM) &&  \text{($D$ is on $X'$ and is exceptional$/X$)},
\end{align*}
so $\mult_DM=0$. Thus $a(D,\Ff,B)=-\epsilon_{\Ff}(D)$. This implies (2). Since $g$ is a divisorial contraction of a prime divisor $E$, by the negativity lemma, $-E$ is ample$/X'$. This implies (3). Since the center of $E$ on $Y$ is a divisor that is exceptional$/X$, $a(E,\Ff,B+M)=-\epsilon_{\Ff}(E)$. Since $E$ is contracted by the $(K_{\Ff_Y}+B_Y+tM_Y+\sum_{i=1}^n\epsilon_{\Ff}(E_i)E_i)$-MMP$/X$: $Y\dashrightarrow X'$, $$a(E,\Ff,B+tM)=a(E,\Ff',B'+tM'+F')>-\epsilon_{\Ff}(E),$$
so $a(E,\Ff,B)>-\epsilon_{\Ff}(E)$. This implies (4).
\end{proof}

\begin{proof}[Proof of Theorem \ref{thm: global acc alg int foliation}]
By Proposition \ref{prop: ACSS modification include some divisor}, possibly replacing $\Ii$ by $\Ii\cup\{1\}$, we may assume that $(X,\Ff,B)$ is $\Qq$-factorial ACSS with a properly associated contraction $\pi: X\rightarrow Z$. Since $(X,\Ff,B)$ is lc, all components of $B$ are horizontal$/Z$. Let $F$ be a general fiber of $\pi$, then $K_{\Ff}|_F=K_X|_F=K_F$ and $(F,B_F:=B|_F)$ is lc. Since $\rk\Ff=r$, $\dim F=r$. By \cite[Theorem 1.5]{HMX14}, the coefficients of $B_F$ belong to a finite set depending only on $r$ and $\Ii$, so the coefficients of $B$ belong to a finite set depending only on $r$ and $\Ii$. The theorem follows.
\end{proof}

\begin{proof}[Proof of Theorem \ref{thm: acc lct alg int foliation}]
Suppose that the theorem does not hold. Then there exists a sequence of lc foliated triples $(X_i,\Ff_i,B_i)$ and $\Rr$-Cartier $\Rr$-divisors $D_i$ on $X_i$, such that $\rk\Ff_i=r$, $B_i,D_i\in\Ii$, and $t_i:=\lct(X_i,\Ff_i,B_i;D_i)$ is strictly increasing. By Lemma \ref{lem: find nontrivial divisor on ACSS model}, possibly replacing $\Ii$ with $\Ii\cup\{1\}$, we may assume that 
\begin{itemize}
    \item[(i)] $(X_i,\Ff_i,B_i+t_i'D_i)$ is $\Qq$-factorial ACSS for some $0<t_i'<t_i$  with a properly associated contraction $\pi_i: X_i\rightarrow Z_i$,
    \item[(ii)] there exists a divisorial contraction $f_i: Y_i\rightarrow X_i$ of a prime divisor $E_i$, such that $\mult_{E_i}D_i>0$ and $a(E_i,\Ff_i,B_i+t_iD_i)=-\epsilon_{\Ff_i}(E_i)$, and
    \item[(iii)] let $B_{Y_i},D_{Y_i}$ be the strict transforms of $B_i,D_i$ on $Y_i$ respectively, then
    $$(Y_i,\Ff_{Y_i}:=f_i^{-1}\Ff_i,B_{Y_i}+t_i'D_{Y_i}+\epsilon_{\Ff_i}(E_i))/Z_i$$
    is $\Qq$-factorial ACSS.
\end{itemize}
Possibly cutting $\Ff_i$ by general elements in base-point-free linear systems as in Propositions \ref{prop: general hyperplane invariant} and \ref{prop: general hyperplane non-invariant}, we may assume that $\Center_{X_i}E_i$ is a closed point. We remark the following: when cutting $\Ff_i$, 
\begin{itemize}
    \item the algebraically integrability of $\Ff_i$ is preserved when cutting $\Ff_i$ by base-point-free linear systems by Proposition \ref{prop: a.i preserved adjunction}, and
    \item (i-iii) are preserved by Propositions \ref{prop: general hyperplane invariant}(3), \ref{prop: general hyperplane non-invariant}(3), and the generality of the hyperplanes.
\end{itemize} 
Let $E_i^\nu$ be the normalization of $E_i$, and
$$K_{\Ff_{E_i}}+B_{E_i}(t):=(K_{\Ff_{Y_i}}+B_{Y_i}+tD_{Y_i}+\epsilon_{\Ff_i}(E_i))|_{E_i^\nu}$$
for any real number $t$. Then $K_{\Ff_{E_i}}+B_{E_i}(t_i)\equiv 0$. By Theorems \ref{thm: adj to invariant} and \ref{thm: adj to non-invariant}, $(E_i^\nu,\Ff_{E_i},B_{E_i}(t_i))$ is lc,
$$K_{\Ff_{E_i}}+B_{E_i}(t)=K_{\Ff_{E_i}}+\sum_{k=1}^{n_i}\frac{w_{i,k}-1+\alpha_{i,k}+t\beta_{i,k}}{w_{i,k}}C_{i,k}+\sum_{k=1}^{l_i}T_{i,k}$$
where $C_{i,k},T_{i,k}$ are prime divisors, $w_{i,k}$ are positive integers, and $\alpha_{i,k},\beta_{i,k}\in\Ii_+$. By Proposition \ref{prop: a.i preserved adjunction}, $\rk\Ff_{E_i}\leq r$. By Theorem \ref{thm: global acc alg int foliation}, possibly passing to a subsequence, we may assume that $\beta_{i,k}=0$ for any $i,k$. Thus $K_{\Ff_{E_i}}+B_{E_i}(t)\equiv 0$ for any $i$ and $t$. Since $f_i$ extracts a unique divisor $E_i$ and $\mult_{E_i}D_i>0$, $K_{\Ff_{Y_i}}+B_{Y_i}+tD_{Y_i}+\epsilon_{\Ff_i}(E_i)E_i$ is ample$/X_i$ when $t>t_i$, so $K_{\Ff_{E_i}}+B_{E_i}(t)$ is ample when $t>t_i$, a contradiction.
\end{proof}

\begin{proof}[Proof of Corollary \ref{cor: global acc rank 1 foliation}]
We may assume that $B\not=0$. By Theorem \ref{thm: subfoliation algebraic integrable}, $\Ff$ is algebraically integrable. Corollary \ref{cor: global acc rank 1 foliation} follows from Theorem \ref{thm: global acc alg int foliation}.
\end{proof}

\begin{proof}[Proof of Theorem \ref{thm: log abundance algebraically integrable foliation}]
By Proposition \ref{prop: ACSS modification include some divisor}, we may assume that $(X,\Ff,B)$ is $\Qq$-factorial ACSS with a properly associated contraction $\pi: X\rightarrow Z$. Let $G$ be a divisor properly associated to $(X,\Ff,B)/Z$. 

By \cite[Proposition 3.6]{ACSS21}, $K_{\Ff}+B\sim_{\mathbb R,Z}K_X+B+G$. Since $K_{\Ff}+B\equiv 0$, $K_X+B+G\equiv_Z0$. Since $(X,B+G)$ is lc, by the theory of uniform rational polytopes (cf. \cite[Lemma 5.3]{HLS19}) and by applying \cite[Corollary 1.6]{HX16}, we have $K_X+B+G\sim_{\mathbb R,Z}0$. Thus $K_{\Ff}+B\sim_{\mathbb R,Z}0$.

By \cite[Theorem 1.3]{LLM23}, we have an lc gfq $(Z,\Ff_Z,B_Z,\Mm)$ induced by a canonical bundle formula $\pi: (X,\Ff,B)\rightarrow Z$. We have
$$K_{\Ff}+B\sim_{\mathbb R}\pi^*(K_{\Ff_Z}+B_Z+\Mm_Z),$$
so $K_{\Ff_Z}+B_Z+\Mm_Z\equiv 0$. By \cite[Theorem 1.3(3)]{LLM23} and \cite[Proposition 3.6(2)]{ACSS21}, there exists a glc g-pair $(Z,B_Z',\Mm)$ induced by a canonical bundle formula $\pi: (X,B+G)\rightarrow Z$. By \cite[Theorem 1.4]{JLX22} (\cite[Theorem 1.1]{FG14} for the $\Qq$-coefficients case), $\Mm$ is $\bb$-abundant. Since $K_{\Ff_Z}=0$, $B_Z=0$ and $\Mm_Z\equiv 0$. Thus $\Mm_Z\sim_{\mathbb R}0$, and $K_{\Ff}+B\sim_{\mathbb R}\pi^*\Mm_Z\sim_{\mathbb R}0$.
\end{proof}

\section{Uniform lc rational polytopes and accumulation points}\label{sec: ulrp}

In this section, we prove Theorem \ref{thm: uniform rational polytope foliation intro} and Corollary \ref{cor: accumulation point ai foliation lct}.

\begin{defn}\label{defn: r affine functional divisor}
Let $X$ be a normal variety, $D_i$ distinct prime divisors, and $d_i(t):\mathbb R\rightarrow\mathbb R$ $\Rr$-affine function. Then we call the formal finite sum $\sum d_i(t)D_i$ an \emph{$\Rr$-affine functional divisor}.
\end{defn}

\begin{defn}
Let $c$ be a non-negative real number, and $\Ii\subset[0,+\infty)$ a set of real numbers. For any $\Rr$-affine functional divisor $\Delta(t)=\sum_{i} d_i(t)D_i$, we write $\Delta(t)\in\mathcal{D}_c(\Ii)$ when the following conditions are satisfied.
\begin{enumerate}
	\item For any $i$, either $d_i(t)=1$, or $d_i(t)$ is of the form $\frac{m-1+\gamma+kt}{m}$, where $m\in\mathbb{N}^{+}$, $\gamma\in\Ii_+$, $k\in\mathbb{Z}$, and
	\item $f+kt$ above can be written as $f+kt=\sum_{j}(f_j+k_jt)$, where $f_j\in\Ii\cup\{0\}$, $k_j\in\mathbb{Z}$, and $f_j+k_jc\ge0$ holds for any $j$.
\end{enumerate}
\end{defn}

\begin{defn}
Let $d$ be a positive integer and  $\Ii\subset[0,+\infty)$ a set of real numbers. We define $\mathcal{B}_{d}(\Ii),\mathcal{B}'_{d}(\Ii)\subset [0,+\infty)$ as follows: $c\in\mathcal{B}_{d}(\Ii)$ (resp. $\mathcal{B}'_{d}(\Ii)$) if and only if there exist a normal projective variety $X$ (resp. a $\Qq$-factorial normal projective variety $X$) and an $\Rr$-affine functional divisor $\Delta(t)$ on $X$ satisfying the following.
\begin{enumerate}
	\item $\dim X\le d$, 
	\item $\Delta(t)\in\mathcal{D}_{c}(\Ii)$,
	\item $(X,\Delta(c))$ is lc,
	\item $K_X+\Delta(c)\equiv0$, and
	\item $K_X+\Delta(c')\not\equiv 0$ for any $c'\neq c$.
\end{enumerate}
\end{defn}

\begin{defn}
Let $r$ be a positive integer and  $\Ii\subset[0,+\infty)$ a set of real numbers. We define $\mathcal{C}_{r}(\Ii),\mathcal{C}'_{r}(\Ii)\subset [0,+\infty)$ as follows: $c\in\mathcal{C}_{r}(\Ii)$ (resp. $c\in\mathcal{C}_{r}'(\Ii)$) if and only if there exist a normal projective variety $X$ (resp. a $\Qq$-factorial normal projective variety $X$), an algebraically integrable foliation $\Ff$ on $X$, and an $\Rr$-affine functional divisor $\Delta(t)$ on $X$ satisfying the following.
\begin{enumerate}
	\item $\rk\Ff\le r$, 
	\item $\Delta(t)\in\mathcal{D}_{c}(\Ii)$,
	\item $(X,\Ff,\Delta(c))$ is lc,
	\item $K_\Ff+\Delta(c)\equiv0$, and
	\item $K_\Ff+\Delta(c')\not\equiv 0$ for any $c'\not=c$.
\end{enumerate}
\end{defn}

\begin{prop}\label{prop: nak special set equal for foliation}
Let $r$ be a positive integer and $\Ii\subset[0,+\infty)$ a set of real numbers. Then $\mathcal{B}_r(\Ii)=\mathcal{C}_r(\Ii)=\mathcal{B}'_r(\Ii)=\mathcal{C}'_r(\Ii)$.
\end{prop}
\begin{proof}
By considering the foliation $\Ff=T_X$, we have $\mathcal{B}_r(\Ii)\subset\mathcal{C}_r(\Ii)$. By the existence of dlt modifications and Proposition \ref{prop: ACSS modification include some divisor}, $\mathcal{B}_r(\Ii)=\mathcal{B}_r'(\Ii)$ and  $\mathcal{C}_r(\Ii)=\mathcal{C}_r'(\Ii)$. We only need show that $\mathcal{C}_r(\Ii)\subset\mathcal{B}_r(\Ii)$.

Pick $c\in\mathcal{C}_r(\Ii)$. The there exists a $\Qq$-factorial normal projective variety $X$, an algebraically integrable foliation $\Ff$ on $X$, and an $\Rr$-affine functional divisor $\Delta(t)$ on $X$, such that
\begin{enumerate}
	\item $\rk\Ff\le r$, 
	\item $\Delta(t)\in\mathcal{D}_{c}(\Ii)$,
	\item $(X,\Ff,\Delta(c))$ is lc,
	\item $K_{\Ff}+\Delta(c)\equiv0$, and
	\item $K_{\Ff}+\Delta(t)\not\equiv 0$ for any $t\not=c$.
\end{enumerate}
By Proposition \ref{prop: ACSS modification include some divisor}, we may let $f: X'\rightarrow X$ be an ACSS modification of $(X,\Ff,\Delta(c))$, $\Ff':=f^{-1}\Ff$, $E:=\Exc(f)^{\Ff'}$, and $\Delta'(t):=f^{-1}_*\Delta(t)+E$ for any real number $t$. Then $\rk\Ff'\leq r$, $\Delta'(t)\in\mathcal{D}_c(\Ii)$, $(X',\Ff',\Delta'(c))$ is lc, and $K_{\Ff'}+\Delta'(c)\equiv 0$. Moreover, for any $t\not=c$, since $$0\not\equiv K_{\Ff}+\Delta(t)=f_*(K_{\Ff'}+\Delta'(t)),$$ $K_{\Ff'}+\Delta'(t)\not\equiv 0$. Therefore, we may replace $(X,\Ff,\Delta(t))$ with $(X',\Ff',\Delta'(t))$, and assume that $(X,\Ff,\Delta(c))$ is $\Qq$-factorial ACSS. Let $\pi: X\rightarrow Z$ be a properly associated contraction of $(X,\Ff,\Delta(c))$ and let $G$ be a divisor properly associated to $(X,\Ff,\Delta(c))/Z$. 

Suppose that for any $0<\delta\ll 1$, $(X,\Ff,\Delta(c+\delta);G)/Z$ or $(X,\Ff,\Delta(c-\delta);G)/Z$ is not ACSS. By Lemmas \ref{lem: acss smaller coefficient} and \ref{lem: acss f-triple perturb coefficient}, there exists a component $D$ of $\Delta(c)$, such that $\mult_D\Delta(c)=1$ and $\mult_D\Delta(t)\not=1$ for any $t\not=c$. By \cite[Lemma 3.7]{Nak16}, $c\in\mathcal{B}_1(\Ii)\subset\mathcal{B}_r(\Ii)$. Therefore, we may assume that $(X,\Ff,\Delta(c+\delta);G)/Z$ and $(X,\Ff,\Delta(c-\delta);G)/Z$ are ACSS for any $0<\delta\ll 1$. 

Fix $0<\delta\ll 1$. Since $K_{\Ff}+\Delta(t)\not\equiv 0$ for any $t\not=c$ and $K_{\Ff}+\Delta(c)\equiv 0$, either $K_{\Ff}+\Delta(c+\delta)$ or $K_{\Ff}+\Delta(c-\delta)$ is not pseudo-effective. We run a $(K_{\Ff}+\Delta(c+\delta))$-MMP (resp. $(K_{\Ff}+\Delta(c-\delta))$-MMP) with scaling of an ample divisor if $(K_{\Ff}+\Delta(c+\delta))$ (resp. $(K_{\Ff}+\Delta(c-\delta))$) is not pseudo-effective. By Lemma \ref{lem: ACSS mmp can run}(4), this MMP is a  $(K_{\Ff}+\Delta(c+\delta)+G)$-MMP$/Z$ (resp. $(K_{\Ff}+\Delta(c-\delta)+G)$-MMP$/Z$) which terminates with a Mori fiber space $\pi'': X''\rightarrow T$ over $Z$. 

We let $\Delta''(t)$ and $G''$ be the images of $\Delta(t)$ and $G$ on $X''$ for any $t$, and let $\Ff''$ be the pushforward of $\Ff$ on $X''$. Then $K_{\Ff''}+\Delta''(c)\equiv 0$ and  $K_{\Ff''}+\Delta''(c+\delta)$ (resp. $K_{\Ff''}+\Delta''(c-\delta)$) is anti-ample$/T$.

Let $F$ be a general fiber of $X''\rightarrow T$. By \cite[Theorem 3.9]{ACSS21}, $F$ is tangent to $\Ff''$, so $K_{\Ff''}|_F=K_{X''}|_F=K_F$. Let $\Delta_{F}(t):=\Delta''(t)|_F$, then
\begin{itemize}
   \item $\dim F\leq r$,
    \item $\Delta_F(t)\in\mathcal{D}_c(\Ii)$,
    \item $(F,\Delta_F(c))$ is lc,
    \item $K_F+\Delta_F(c)\equiv 0$, and
    \item $K_F+\Delta_F(c+\delta)$ or $K_F+\Delta_F(c-\delta)$ is anti-ample.
\end{itemize}
Thus $c\in\mathcal{B}_r(\Ii)$.
\end{proof}

\begin{thm}\label{thm: uniform rational polytope foliation one variable}
Let $d,c,m$ be positive integers, $r_1,\dots,r_c$ real numbers such that $1,r_1,\dots,r_c$ are linearly independent over $\mathbb Q$, $\bm{r}:=(r_1,\dots,r_c)$, and $s_1,\dots,s_m: \mathbb R^{c+1}\rightarrow\mathbb R$ $\mathbb Q$-linear functions. Then there exists a positive real number $\delta$ depending only on $\bm{r}$ and $s_1,\dots,s_m$ satisfying the following. Assume that
\begin{enumerate}
    \item $(X,\Ff,B=\sum_{i=1}^ms_i(1,r_1,\dots,r_{c-1},t)B_i)$ is an lc f-triple such that $\Ff$ is algebraically integrable and $\rk\Ff\leq d$,
    \item $B_i\geq 0$ are distinct Weil divisors (possibly $0$) and $s_i(1,\bm{r})\geq 0$, and
    \item $B(t):=\sum_{i=1}^ms_i(1,r_1,\dots,r_{c-1},t)B_i$ for any $t\in\mathbb R$.
\end{enumerate}
Then $(X,\Ff,B(t))$ is lc for any $t\in (r_c-\delta,r_c+\delta)$.
\end{thm}
\begin{proof}
We let $s_i(t):=s_i(1,r_1,\dots,r_{c-1},t)$ for any $t\in\mathbb R$. If $s_i(r_c)=0$, then $s_i(t)=0$ for any $i$, so we may assume that $s_i(r_c)\not=0$ for any $i$. By Proposition \ref{prop: ACSS modification include some divisor}, possibly replacing $(X,\Ff,B(r_c))$ with an ACSS model, we may assume that $(X,\Ff,B(r_c))$ is $\Qq$-factorial ACSS.

We only need to prove that there exists a positive real number $\epsilon$ depending only on $\Ii$ and $r$, such that for any lc threshold $t_0$ of $(X,\Ff,B(t))$, $|t_0-r_c|>\epsilon$. Thus we may assume that $(X,\Ff,B(t))$ has an lc threshold $t_0$. Since $1,r_1,\dots,r_c$ are linearly independent over $\Qq$, $r_c\not=t_0$. Moreover, there exists a positive real number $\delta_1$, such that for any $t\in (r_c-\delta_1,r_c+\delta_1)$ and any $i$, $s_i(t)\geq\frac{1}{2}s_i(r_c)>0$. In particular, for any $t\in (r_c-\delta_1,r_c+\delta_1)$, $\Supp B(t)=\Supp B(r_c)$, $\Supp\lfloor B(t)\rfloor=\Supp\lfloor B(r_c)\rfloor$, and $B(t)\geq\frac{1}{2}B(r_c)$.

We may assume that $t_0\in (r_c-\delta_1,r_c+\delta_1)$. In particular, we may assume that $(X,\Ff,B(t_0))$ does not have an lc center in codimension $1$. Thus $(X,\Ff,B(t_0))$ has an lc center $x$ such that $\dim x\leq \dim X-2$, and $x$ is not an lc center of $(X,\Ff,B(r_c))$. In particular, $x\in\Supp B(r_c)$.

By Lemma \ref{lem: find nontrivial divisor on ACSS model}, possibly replacing $(X,\Ff,B(t))$, we may assume that there exist a divisorial contraction $g: Y\rightarrow X$ of a prime divisor $E$ and a real number $s$ satisfying the following: let $B_Y(t)$ be the strict transform of $B(t)$ on $Y$ for any $t$ and $\Ff_Y:=g^{-1}\Ff$, then
\begin{itemize}
\item[(i)] $s\in (r_c,t_0)$ if $r_c>t_0$, and $s\in (t_0,r_c)$ if $t_0<r_c$,
\item[(ii)] $(X,\Ff,B(s))$ is $\Qq$-factorial ACSS, $(X,\Ff,B(r_c))$ is lc, and  $(X,\Ff,B(t_0))$ is lc,
\item[(iii)] $-E$ is ample over $X$,
\item[(iv)] $(Y,\Ff_Y,B_Y(s)+\epsilon_{\Ff}(E))$ is $\Qq$-factorial ACSS, and
\item[(v)] $a(E,\Ff,B(t_0))=-\epsilon_{\Ff}(E)$ and $a(E,\Ff,B(r_c))>-\epsilon_{\Ff}(E)$. In particular,  $(Y,\Ff_Y,B_Y(t_0)+\epsilon_{\Ff}(E))$ is lc.
\end{itemize}
Possibly cutting $\Ff$ by general elements in base-point-free linear systems as in Propositions \ref{prop: general hyperplane invariant} and \ref{prop: general hyperplane non-invariant}, we may assume that $\Center_XE$ is a closed point. We remark the following: when cutting $\Ff_i$, 
\begin{itemize}
    \item the algebraically integrability of $\Ff$ is preserved when cutting $\Ff$ by sections by Proposition \ref{prop: a.i preserved adjunction}, and
    \item (i-v) are preserved by Propositions \ref{prop: general hyperplane invariant}(3), \ref{prop: general hyperplane non-invariant}(3), and the generality of the hyperplanes.
\end{itemize} 
Let $E^\nu$ be the normalization of $E$, and
$$K_{\Ff_{E}}+B_{E}(t):=(K_{\Ff_{Y}}+B_Y(t)+\epsilon_{\Ff}(E))|_{E^\nu}$$
for any real number $t$. Then $K_{\Ff_{E}}+B_{E}(t_0)\equiv 0$ and $K_{\Ff_{E}}+B_{E}(t)\not\equiv 0$ for any $t\not=t_0$. By Theorems \ref{thm: adj to invariant} and \ref{thm: adj to non-invariant}, $(E^\nu,\Ff_{E},B_{E}(t_0))$ is lc and $B_E(t)\in\mathcal{D}_{t_0}(\Ii)$. By Proposition \ref{prop: a.i preserved adjunction}, $\Ff_E$ is algebraically integrable and $\rk\Ff_E\leq d$. By Proposition \ref{prop: nak special set equal for foliation}, $t_0\in\mathcal{B}'_d(\Ii)$. The theorem follows from \cite[Corollary 3.9]{Nak16}.
\end{proof}

\begin{thm}\label{thm: uniform rational polytope}
Let $d,c,m$ be positive integers, $r_1,\dots,r_c$ real numbers such that $1,r_1,\dots,r_c$ are linearly independent over $\mathbb Q$, $\bm{r}:=(r_1,\dots,r_c)$, and $s_1,\dots,s_m: \mathbb R^{c+1}\rightarrow\mathbb R$ $\mathbb Q$-linear functions. Then there exists an open subset $U\ni\bm{r}$ depending only on $\bm{r}$ and $s_1,\dots,s_m$ satisfying the following. Assume that
\begin{enumerate}
    \item $(X,\Ff,B(\bm{r}):=\sum_{i=1}^ms_i(1,\bm{r})B_i)$ is an lc f-triple such that $\Ff$ is algebraically integrable and $\rk\Ff\leq d$,
    \item $B_i\geq 0$ are distinct Weil divisors (possibly $0$) and $s_i(1,\bm{r})\geq 0$, and
    \item $B(\bm{v}):=\sum_{i=1}^ms_i(1,\bm{v})B_i$ for any $t\in\mathbb R$.
\end{enumerate}
Then $(X,\Ff,B(\bm{v}))$ is lc for any $\bm{v}\in U$.
\end{thm}
\begin{proof}
We apply induction on $c$. When $c=1$, Theorem \ref{thm: uniform rational polytope} directly follows from Theorem \ref{thm: uniform rational polytope foliation one variable}. When $c\geq 2$, by Theorem \ref{thm: uniform rational polytope foliation one variable}, there exists a positive integer $\delta$ depending only on $r_1,\dots,r_c,s_1,\dots,s_m$, such that for any $t\in (r_c-\delta,r_c+\delta)$, $(X,\Ff,\sum_{i=1}^ms_i(1,r_1,\dots,r_{c-1},t)B_i)$ is lc. We pick rational numbers $r_{c,1}\in (r_c-\delta,r_c)$ and $r_{c,2}\in (r_c,r_c+\delta)$ depending only on $r_1,\dots,r_c,s_1,\dots,s_m$. By induction on $c$, there exists an open subset $U_0\in (r_1,\dots,r_{c-1})$ of $\mathbb R^{c-1}$, such that for any $\bm{v}\in U_0$, $(X,\Ff,\sum_{i=1}^ms_i(1,\bm{v},r_{c,1})B_i)$ and $(X,\Ff,\sum_{i=1}^ms_i(1,\bm{v},r_{c,2})B_i)$ are lc. We may pick $U:=U_0\times (r_{c,1},r_{c,2})$.
\end{proof}

\begin{proof}[Proof of Theorem \ref{thm: uniform rational polytope foliation intro}]
It immediately follows from Theorem \ref{thm: uniform rational polytope}.
\end{proof}

\begin{thm}\label{thm: accumulation point of foliated lc threshold complicated version}
Let $d,c$ be two non-negative integers, $r_1,\dots,r_c$ real numbers, and $\Ii\subset [0,1]$ a DCC set, such that $\bar\Ii\subset\Span_{\mathbb Q}(1,r_1,r_2,\dots,r_c)$. The the accumulation points of 
$$\{\lct(X,\Ff,B;D)\mid \rk\Ff\leq d,\Ff\text{ is algebraically integrable}, (X,\Ff,B)\text{ is lc, }B\in\Ii,D\in\mathbb N^+\}$$
belong to $\Span_{\mathbb Q}(1,r_1,r_2,\dots,r_c)$.
\end{thm}
\begin{proof}
Suppose the theorem does not hold. By Theorem \ref{thm: acc lct alg int foliation}, there exist a sequence of lc f-triples $(X_i,\Ff_i,B_i)$ and effective $\Qq$-Cartier Weil divisors $D_i$ on $X_i$, such that $\Ff_i$ is algebraically integrable, $\rk\Ff_i\leq r$, $t_i:=\lct(X_i,\Ff_i,B_i;D_i)$ is strictly decreasing, $B_i\in\Ii$, and $t:=\lim_{i\rightarrow+\infty}t_i\not\in\Span_{\mathbb Q}(1,r_1,r_2,\dots,r_c)$. We write $B_i=\sum_{j=1}^{m_i}b_{i,j}B_{i,j}$, where $B_{i,j}$ are the irreducible components of $B_i$. Possibly replacing $(X_i,\Ff_i,B_i)$ with an ACSS model and replacing $D_i$ with its pullback, we may assume that $(X_i,\Ff_i,B_i)$ is $\Qq$-factorial ACSS. 

We let $G_i$ be a divisor properly associated to $(X_i,\Ff_i,B_i)$, then $(X_i,B_i+G_i)$ is lc. We let $x_i\in X_i$ be a point such that $(X_i,\Ff_i,B_i+t_iD_i)$ is lc near $x_i$ and $(X_i,\Ff_i,B_i+(t_i+\delta)D_i)$ is not lc near $x_i$ for any positive real number $\delta$. Since  $(X_i,B_i+G_i)$ is lc near $x_i$, by \cite[18.22 Theorem]{Kol+92}, possibly shrinking $X_i$ to a neighborhood of $x_i$ and passing to a subsequence, $m_i$ is a constant for every $i$ an we denote it by $m$. 

Possibly passing to a subsequence again, we may assume that $b_{i,j}$ is strictly increasing for every fixed $j$. We let $\bar b_j:=\lim_{i\rightarrow+\infty}b_{i,j}$ and $\bar B_i:=\sum_{j=1}^m\bar b_jB_{i,j}$.

Since $(X_i,B_i+t_iD_i)$ is lc, $(X_i,B_i+tD_i)$ is lc. By Theorem \ref{thm: acc lct alg int foliation}, possibly passing to a subsequence, we may assume that $(X_i,\bar B_i+tD_i)$ is lc for each $i$. By our assumption, $\bar b_j\in\Span_{\mathbb Q}(1,r_1,\dots,r_c)$ for each $j$ and $t\not\in\Span_{\mathbb Q}(1,r_1,\dots,r_c)$. We let $V$ be the rational envelope of $(\bar b_1,\dots,\bar b_m,t)$ in $\mathbb R^{m+1}$, then $V=V'\times\mathbb R$, where $V'$ is the rational envelope of $(\bar b_1,\dots,\bar b_m)$ in $\mathbb R^m$. By Theorem \ref{thm: uniform rational polytope foliation intro}, there exist an open subset $U'\ni (\bar b_1,\dots,\bar b_m)$ of $V'$, and an open subset $W\ni t$ of $\mathbb R$, such that $(X_i,\sum_{j=1}^mv_jB_{i,j}+wD_i)$ is lc for any $(v_1,\dots,v_m)\in U'$ and $w\in W$. In particular, $(X_i,\bar B_i+wD_i)$ is lc for any $w\in W$. Possibly passing to a subsequence, we may assume that there exists a real number $w_0\in W$ such that $w_0>t_i$ for any $i$. Then $(X_i,\bar B_i+w_0D_i)$ is lc for any $i$, so $(X_i,B_i+w_0D_i)$ is lc for any $i$, so
$$t_i=\lct(X_i,\Ff_i,B_i;D_i)\geq w_0>t_i,$$
a contradiction.
\end{proof}

\begin{proof}[Proof of Corollary \ref{cor: accumulation point ai foliation lct}]
It immediately follows from Theorem \ref{thm: accumulation point of foliated lc threshold complicated version} by taking $c=0$.
\end{proof}

\end{document}